\documentclass[a4paper,12pt]{article}
\RequirePackage{amsthm,amsmath,amsfonts,amssymb}
\RequirePackage[numbers]{natbib}
\RequirePackage{graphicx}
\usepackage{authblk}
\usepackage[dvipsnames]{xcolor}
\usepackage{latexsym,dsfont,enumerate,bbm,threeparttable,rotating, pdflscape,verbatim,amssymb}
\usepackage{algorithm,algcompatible,algpseudocode,booktabs,xcolor,tikz}
\usetikzlibrary{arrows,shapes,chains}
\usepackage[normalem]{ulem}

\usepackage[colorlinks,linkcolor=red,citecolor=blue,urlcolor=blue,breaklinks]{hyperref}
\usepackage{cleveref}

\usepackage{threeparttable}
\newcommand{\bQ}{\mbox{\bf Q}}

\newcommand{\bX}{\mbox{\bf X}}
\newcommand{\bH}{\mbox{\bf H}}
\newcommand{\bP}{\mbox{\bf P}}
\newcommand{\bI}{\mbox{\bf I}}
\newcommand{\by}{\mbox{\bf y}}
\newcommand{\bx}{\mbox{\bf x}}

\newcommand{\bz}{\mbox{\bf z}}
\newcommand{\bt}{\mbox{\bf t}}

\newcommand{\bb}{\mbox{\bf b}}

\newcommand{\bv}{\mbox{\bf v}}
\newcommand{\bu}{\mbox{\bf u}}

\newcommand{\bkappa}{\mbox{\boldmath $\kappa$}}
\newcommand{\bh}{\mbox{\bf h}}

\newcommand{\bzero}{\mbox{\bf 0}}
\newcommand{\bPhi}{\mbox{\boldmath $\Phi$}}
\newcommand{\bSigma}{\mbox{\boldmath $\Sigma$}}
\newcommand{\bphi}{\mbox{\boldmath $\phi$}}
\newcommand{\bveps}{\mbox{\boldmath $\varepsilon$}}

\newcommand{\bbeta}{\mbox{\boldmath $\beta$}}
\newcommand{\bpsi}{\mbox{\boldmath $\psi$}}
\newcommand{\bomega}{\mbox{\boldmath $\omega$}}
\newcommand{\bgamma}{\mbox{\boldmath $\gamma$}}

\newcommand{\brho}{\mbox{\boldmath $\rho$}}
\newcommand{\bTheta}{\mbox{\boldmath $\Theta$}}

\newcommand{\hbeta}{\widehat{\beta}}
\newcommand{\bpi}{\mbox{\boldmath $\pi$}}

\newcommand{\homega}{\widehat{\omega}}
\newcommand{\sgn}{\mathrm{sgn}}

\newtheorem{theorem}{Theorem}
\newtheorem{lemma}{Lemma}
\theoremstyle{remark}
\newtheorem{corollary}{Corollary}
\newtheorem{definition}{Definition}
\newtheorem{proposition}{Proposition}
\newtheorem{condition}{Condition}

\providecommand{\keywords}[1]{\textbf{\textit{Keywords---}} #1}

\begin{document}

%

\title{High-dimensional inference via hybrid orthogonalization
	\thanks{Yang Li is Postdoctoral Fellow, International Institute of Finance, The School of Management, University of Science and Technology of China, Hefei, Anhui, 230026, China (E-mail:tjly@mail.ustc.edu.cn). Zemin Zheng is Professor, International Institute of Finance, The School of Management, University of Science and Technology of China, Hefei, Anhui, 230026, China (E-mail:zhengzm@ustc.edu.cn). Jia Zhou is Ph.D. candidate, International Institute of Finance, The School of Management, University of Science and Technology of China, Hefei, Anhui, 230026, China (E-mail:tszhjia@mail.ustc.edu.cn).
Ziwei Zhu is Assistant Professor, Department of Statistics, University of Michigan, Ann Arbor, MI 48109 (E-mail:ziweiz@umich.edu).
}}
\author[1]{Yang Li}
\author[1]{Zemin Zheng}
\author[1]{Jia Zhou}
\author[2]{Ziwei Zhu}

\affil[1]{University of Science and Technology of China}
\affil[2]{University of Michigan, Ann Arbor}
\date{November 26, 2021}
\maketitle

\begin{abstract}
	The past decade has witnessed a surge of endeavors in statistical inference for high-dimensional sparse regression, particularly via de-biasing or relaxed orthogonalization. Nevertheless, these techniques typically require a more stringent sparsity condition than needed for estimation consistency, which seriously limits their practical applicability. To alleviate such constraint, we propose to exploit the identifiable features to residualize the design matrix before performing debiasing-based inference over the parameters of interest. This leads to a hybrid orthogonalization (HOT) technique that performs strict orthogonalization against the identifiable features but relaxed orthogonalization against the others. Under an approximately sparse model with a mixture of identifiable and unidentifiable signals, we establish the asymptotic normality of the HOT test statistic while accommodating as many identifiable signals as consistent estimation allows. The efficacy of the proposed test is also demonstrated through simulation and analysis of a stock market dataset.
\end{abstract}


\section{Introduction}

Feature selection is crucial to a myriad of modern machine learning problems under high dimensions, including bi-clustering analysis with high-throughput genomic data, collaborative filtering in sophisticated recommendation systems, compression of deep neural networks, among others. However, characterizing the uncertainty of feature selection is challenging: the curse of high dimensionality of covariates can undermine the validity of classical inference procedures  \citep{Fan19, Sur20}, which only holds under low-dimensional setups.

When the dependence structure of the response variable on the explanatory variables takes some special form such as a linear one, a series of statistical inference methods have been proposed to mitigate the curse of dimensionality and derive valid p-values for these explanatory variables.
Among them, a line of research performs statistical inference after model selection with significance statements conditional on the selected model \citep{Lockhart14, Chernozhukov15, Lee16, Tibshirani2016, Tian18}. Another complementary class of methods, which are more closely related with this work, propose to de-bias the penalized estimators by inverting the Karush-Kuhn-Tucker condition of the corresponding optimization problem or equivalently through low-dimensional projections \citep{Javanmard2014, vandeGeer2014, Zhang2014}, to establish unconditional significance statements. The major advantage of such inference procedures is that they do not require the true signals to be uniformly strong, thereby accommodating general magnitude of true regression coefficients. This de-biasing technique has been extended to inference of graphical structures \citep{Jankova2015, RenZhou15, Chen16, Jankova2017, Ren19}, goodness-of-fit tests \citep{Jankova2020} and empirical likelihood tests \citep{Chang20}. It has also been generalized to remove the impact of high-dimensional nuisance parameters \cite{Chernozhukov18} and tackle distributed high-dimensional statistical inference \citep{Lee17, Zhu2018, Lian18, Zheng2020}. 

There is no free lunch though: the asymptotic normality of the de-biased estimators comes at the price of a strict sparsity constraint: $s = o(\sqrt n / \log p)$, where $n$ is the sample size, $p$ is the dimension of the parameter of interest, and $s$ denotes the number of nonzero coefficients. This is more restrictive than $s = o(n / \log p)$, which suffices for consistent estimation in $L_2$-loss.
When the true regression coefficients are strictly sparse ($L_0$-sparse) and the row-wise maximum sparsity $s_0$ of the population precision matrix of the features satisfies $s_0 = o(\sqrt{n}/ \log p)$,  \cite{Javanmard2018} developed a novel leave-one-out proof technique to show that the constraint on $s$ can be relaxed to a nearly optimal one: $s = o\{n/ (\log p)^2\}$. For approximately sparse (capped $L_1$-sparse, see \eqref{eq:spars}) structures, it remains unknown if the sparsity constraint can be relaxed similarly.

Our work is mainly motivated by a key observation that despite of embracing signals of various strength, the existing de-biasing inference methods \citep{Zhang2014, Javanmard2014, vandeGeer2014} treat all the covariates equally, regardless of their signal strength. Then the estimation errors on all the signal coordinates can accumulate and jeopardize the asymptotic validity of the resulting tests. Note that when the magnitude of a regression coefficient is of order $\Omega\{\sqrt{(\log p) /n}\}$, it can be identified by feature screening techniques \citep{Fan2008, Li2012, Wang2016} or variable selection methods \citep{Wainwright09, Zhang10, Zhang2014} with high probability.
Then a natural question arises: if those strong signals can be identified in the first place, is the number of them still limited to $o(\sqrt{n}/ \log p)$ for valid inference? If not, what will be the new requirement for the number of identifiable signals, particularly in the presence of many weak signals? Another interesting question is whether the test efficiency can be improved by exploiting these identifiable signals. If so, then how?

In this paper, we intend to address these questions by developing and analyzing a new inference methodology that takes full advantage of the identifiable signals. We propose to first identify the strong signals and eliminate their impact on the subsequent test.
Specifically, we propose a new residualization technique called Hybrid OrThogonalization (HOT) that residualizes the features of our interest such that they are strictly orthogonal to the identifiable features and approximately orthogonal to the remaining unidentifiable ones. The resulting HOT test thus avoids the bias of estimating identifiable coefficients in the original de-biasing step and substantially relaxes the requirement on the sparsity $s$. To the best of our knowledge, such a hybrid inference procedure is new. 


\subsection{Contributions and organization}

We summarize our main contributions of this work as follows:
\begin{enumerate}
	\item We construct the HOT test statistic via exact orthogonalization against strong signals followed by relaxed orthogonalization against the remaining signals. Denote by $s_0$ the maximum row-wise sparsity level of the precision matrix of the predictors. We establish the asymptotic normality of the HOT test statistic, which allows the number of identifiable coefficients $s_1$ to be of order $o(n/ s_0)$. One crucial ingredient of our analysis is verifying that the sign-restricted cone invertibility factor of the residualized design matrix that is yielded by the exact orthogonalization step is well above zero.
	
	\item We propose another equivalent formulation of the HOT test statistic through a partially penalized least squares problem (see Proposition \ref{pro1}). Analyzing this form also leads to the asymptotic normality of the HOT test statistic, but with a different requirement on $s_1$: $s_1 = o(n / \log p)$. This matches the sparsity requirement for consistent estimation in terms of the Euclidean norm. Combining this sparsity constraint with the first one gives a further relaxed constraint: $s_1 = o\{\max(n / \log p, n / s_0)\}$.
	
	\item We show that the HOT test allows the number of weak signals (non-identifiable) to be of order $o(\sqrt{n}/\log p)$, which is the same order as the sparsity constraint under existing de-biasing works.
	
	\item Last but not least, we show that the asymptotic efficiency of the HOT test is the same as the existing de-biasing methods, indicating no efficiency loss due to the hybrid orthogonalization.
\end{enumerate}

The rest of this paper is organized as follows. Section \ref{me} presents the motivation and the new methodology of the HOT-based inference. We establish comprehensive theoretical properties to ensure the validity of the HOT tests in Section \ref{tr}. Simulated data analyses are provided in Section \ref{sr}. We also demonstrate the usefulness of the HOT test by analyzing the data set of stock short interest and aggregate return in Section \ref{rd}. Section \ref{dis} concludes with discussion and future research directions. The proofs and additional technical details are relegated to the supplementary materials to this paper.

\subsection{Notation}
For any vector $\bx=(x_1,\ldots,x_p)^\top$, denote by $\|\bx\|_q=(\sum_{j=1}^p|x_j|^q)^{1/q}$ the $L_q$-norm for $q \in (0, \infty)$.
Given two sets $S_1$ and $S_2$, we use $S_1 \backslash S_2$ to denote $S_1 \cap S_2 ^ c$.
For any matrix $\bX$,   denote by $\bX_{S_1}$  the submatrix of $\bX$ with columns in $S_1$,  and $\bX_{S_1S_2}$  the submatrix of $\bX$ with rows in $S_1$ and columns in  $S_2$. We write $\bX_{-j}$ as the submatrix of $\bX$ with columns constrained on the set $\{1,2,\dots,p\} \backslash \{j\}$.
For two series $\{a_n\}$ and $\{b_n\}$, we write $a_n = O(b_n)$ if there exists a universal constant $C$ such that $a_n \le Cb_n$ when $n$ is sufficiently large. We write $a_n = \Omega(b_n)$ if there exists a universal constant $c$ such that $a_n \ge cb_n$ when $n$ is sufficiently large.
We say $a_n \asymp b_n$ if $a_n = O(b_n)$ and $a_n = \Omega(b_n)$.
\section{Efficient inference via hybrid orthogonalization}\label{me} 

\subsection{Setup and motivation} \label{mot} 

Consider the following high-dimensional linear regression model: 
\begin{equation}\label{eq002p}
\by = \bX\bbeta + \bveps,
\end{equation}
where $\by = (y_1, \ldots, y_n)^\top$  is an $n$-dimensional vector of responses, $\bX = (\bx_1, \ldots, \bx_p)$ is an $n\times p$ design matrix with $p$ covariates, $\bbeta = (\beta_1,\ldots,\beta_p)^\top$ is a $p$-dimensional unknown regression coefficient vector, and $\bveps=(\varepsilon_1,\ldots,\varepsilon_n)^\top$ is an $n$-dimensional random error vector that is independent of $\bX$.
Under a typical high-dimensional sparse setup, the number of features $p$ is allowed to grow exponentially fast with the sample size $n$, while the true regression coefficient vector $\bbeta$ maintains to be parsimonious.  
Given that the strict sparsity and uniform signal strength assumptions are often violated in practice, we impose rather the weaker capped $L_1$ sparsity condition that was suggested in \cite{Zhang2014}:
\begin{equation}
	\label{eq:spars}
    \sum_{j=1}^{p}\min{\big\{|\beta_{j}|/(\sigma \lambda_0),1\big\}}\leq s, 
\end{equation}
where $\lambda_0 = \sqrt{2 \log (p)/n}$ is in accordance with the level of the maximum spurious correlation $\|n^{-1} \bX^{\top} \bveps\|_{\infty}$. Specifically, a signal gets fully counted towards the sparsity allowance only when its order of magnitude exceeds $\sigma\lambda_0$, which allows it to stand out from spurious predictors and be identified through feature screening or selection \citep{Fan2008, Wainwright09, Zhang10}.
The weaker signals instead have discounted contribution to $s$, whose magnitude depends on their strength. Therefore, depending on whether the signal strength exceeds $\sigma\lambda_0$, the sparsity parameter $s$ can naturally be decomposed into two parts $s_1$ and $s_2$ that correspond to the identifiable and weak signals respectively. As we shall see in the sequel, the major bottleneck of our HOT test lies in only $s_2$, while the constraint on $s_1$ by standard de-biasing methods can be substantially relaxed through the proposed hybrid orthogonalization. 


To motivate our approach, we first introduce the low-dimensional projection estimator (LDPE) $\widehat{\bbeta}^{\rm L} = (\hbeta_1^{\rm L}, \cdots, \hbeta_p^{\rm L})^\top$ proposed in \cite{Zhang2014}. Consider an initial estimator $\widehat{\bbeta}^{\text{(init)}}=(\hbeta_1^{\text{(init)}},\ldots,\hbeta_p^{\text{(init)}})^\top$ obtained through the scaled lasso regression \citep{Sun2012} with tuning parameter $\lambda_0$:
\begin{equation}\label{sl}
	(\widehat{\bbeta}^{\text{(init)}}, \widehat{\sigma}) \in  \mathop{\arg\min}_{\bbeta, \sigma}\left\{(2\sigma n)^{-1}\|\by-\bX\bbeta\|_2^2+ 2^{-1}\sigma + \lambda_0\|\bbeta\|_1 \right\}.
\end{equation}
Then for $1 \leq j \leq p$, define the LDPE as
\begin{equation*}
\hbeta_j^{\rm L} := \hbeta_j^{\text{(init)}} + \bz_j^\top (\by - \bX \widehat{\bbeta}^{\text{(init)}})/\bz_j^\top \bx_j,
\end{equation*}
where $\bz_j$ is a relaxed orthogonalization of $\bx_j$ against $\bX_{-j}$ and can be generated as the residual of the corresponding lasso regression \citep{Tibshirani1996}.
To see the asymptotic distribution of LDPE, note that
\begin{eqnarray} \label{jugaoshab}
\hbeta_j^{\rm L} - \beta_j = \frac{\bz_j^\top\bveps}{\bz_j^\top\bx_j} + \sum_{k\ne j}\frac{\bz_j^\top\bx_k(\beta_k-\widehat{\beta}_k^{(\rm init)})}{\bz_j^\top\bx_j}.
\end{eqnarray}
The first term on the right hand side (RHS) dictates the asymptotic distribution of the LDPE and is of order $O_{\mathbb{P}}(n^{-1/2})$.
The second term on the RHS is responsible for the bias of the LDPE. Define the bias factor
\[
    \eta_j := \max_{k\ne j} |\bz_j^\top \bx_k|/\|\bz_j\|_2
\]
and the noise factor
\[
    \tau_j := \|\bz_j\|_2/|\bz_j^\top \bx_j|.
\]
Then the bias term of the LDPE is bounded as
\begin{eqnarray}\label{sbc}
\sum_{k\ne j}\frac{|\bz_j^\top\bx_k(\beta_k-\widehat{\beta}_k^{(\rm init)})|}{|\bz_j^\top\bx_j|} \leq \bigg(\max_{k\ne j} \frac{|\bz_j^\top \bx_k|}{\|\bz_j\|_2}\bigg)  \frac{\|\bz_j\|_2}{|\bz_j^\top \bx_j|} \big\|\widehat{\bbeta}^{\text{(init)}} - \bbeta\big\|_1 = \eta_j \tau_j \|\widehat{\bbeta}^{\text{(init)}} - \bbeta\|_1.
\end{eqnarray}
When the design is sub-Gaussian, $\eta_j \asymp \sqrt{\log p}$ and $\tau_j \asymp n ^ {- 1 / 2}$ under typical realizations of $\bX$. 
Therefore, in order to achieve asymptotic normality of $\widehat\beta_j ^ \mathrm{L}$, we require the bias term in \eqref{sbc} to be $o_{\mathbb{P}}(n^{-1/2})$, which entails that
\begin{eqnarray*}
\eta_j \|\widehat{\bbeta}^{\text{(init)}} - \bbeta\|_1 = o_{\mathbb{P}}(1).
\end{eqnarray*}
Under mild regularity conditions, \cite{Zhang2014} showed that $\|\widehat{\bbeta}^{\text{(init)}} - \bbeta\|_1 = O_{\mathbb{P}}(s \sqrt{\log (p) /n})$. Combining this with the aforementioned order of $\eta_j$ suggests that one  requires
\begin{equation}
    \label{eq:ldpe_s_bound}
    s = o(\sqrt{n} / \log p)
\end{equation}
to achieve the asymptotic normality of $\widehat\beta_j ^{\mathrm{L}}$.
Similar sparsity constraint is necessitated by other de-biasing methods such as \cite{Javanmard2014} and \cite{vandeGeer2014}, implying that valid inference via the de-biased estimators typically requires the true model to be more sparse than that with $L_2$ estimation consistency, which requires merely that $s = o(n / \log p)$. 

The advantage of this inference procedure is that it applies to general magnitudes of the true regression coefficients $\beta_j$ under sparse or approximately sparse structures as long as the aforementioned accuracy of the initial estimate holds. The downside, however, is that it requires the stringent sparsity assumption \eqref{eq:ldpe_s_bound}. Looking back at \eqref{jugaoshab} and the definition of $\bz_j$, we find that the construction of $\bz_j$ in LDPE does not discriminate between identifiable and weak signals. Then even with an initial estimate to offset the impacts of identifiable signals in the bias term, the estimation errors of the $s_1$ identifiable nonzero coefficients can aggregate to be of order $s_1 \sqrt{\log (p) /n}$, which gives rise to a harsh constraint on the number of identifiable predictors. 

To address this issue, we propose to construct the vector $\bz_j$ by a hybrid orthogonalization technique after identifying large coefficients by feature screening or selection procedures. A partial list of such procedures include SIS \citep{Fan2008}, distance correlation learning \citep{Li2012}, HOLP \citep{Wang2016}, lasso \citep{Tibshirani1996}, SCAD \citep{FanLi2001}, the adaptive lasso \citep{Zou2006}, the Dantzig selector \citep{Candes2007}, MCP \citep{Zhang10}, the scaled lasso \citep{Sun2012}, the combined $L_1$ and concave regularization \citep{FanLv2014}, the thresholded regression \citep{Zheng2014}, the constrained Dantzig selector \citep{Kong16}, among many others \citep{James2008}.
Intuitively, an ideal hybrid orthogonalization vector $\bz_j$ should be orthogonal to the identifiable features so that their impact can be completely removed in the bias term even without resorting to any initial estimate. At the same time, it should be a relaxed orthogonalization against the remaining covariate vectors in case that some weak signals are above the magnitude of $n^{-1/2}$ and break down the asymptotic normality.

\subsection{Unveiling the mystery of the hybrid orthogonalization}\label{ms} 

The first step of hybrid orthogonalization is residualizing the target covariate and the unidentifiable ones using the identifiable covariates through exact orthogonalization. After obtaining those residualized covariate vectors, we then construct the relaxed orthogonalization vector via penalized regression. Specifically, given a pre-screened set $\widehat{S}$ of identifiable signals and any target index $1 \le j \le p$, the HOT vector $\bz_j$ can be constructed in two steps: 

\emph{Step 1. Exact orthogonalization.} For $k \in \{j\} \cup \widehat{S}^c$, we take $\bpsi_k^{(j)}$ as the exact orthogonalization of $\bx_k$ against $\bX_{\widehat{S} \backslash \{j\}}$. That is, 
\begin{equation}\label{zdef}
\bpsi_k^{(j)} = (\bI_n-\mathbb{\bP}_{\widehat{S}\backslash \{j\}})\bx_k,
\end{equation}
where $\mathbb{\bP}_A=\bX_A(\bX_A^\top \bX_A)^{-1}\bX_A^\top$ is the projection matrix of the column space of $\bX_A$. Therefore, for any $k \in \{j\} \cup \widehat{S}^c$, the projected vector $\bpsi_k^{(j)}$ is orthogonal to $\bX_{\widehat{S}\backslash  \{j\}}$. The submatrix of the identifiable covariates $\bX_{\widehat{S}\backslash  \{j\}}$ will then be discarded in the construction of $\bz_j$ after this step.
It is worth noticing that when $j \in \widehat{S}^c$, for any $k \in \widehat{S}^c$, $\bpsi_k^{(j)} = (\bI_n - \bP_{\hat S})\bx_k$, which is the residual vector of projecting $\bx_k$ onto the space spanned by all covariate vectors in $\widehat{S}$. However, $\bpsi_k^{(j)}$ varies across $j \in \widehat{S}$, since the corresponding projection space only includes the other covariate vectors in $\widehat{S}$.

\smallskip

\emph{Step 2. Relaxed orthogonalization.} The hybrid orthogonalization vector $\bz_j$ is constructed as the residual of the following lasso regression based on $\{\psi ^ {(j)}_k\}_{k \in \{j\} \cup \widehat S ^ c}$: 
\begin{align}
&\bz_j = \bpsi_j^{(j)} - \bpsi_{\widehat{S}^c\backslash \{j\}}^{(j)} \widehat{\bomega}_j; \label{penl}\\
&\widehat{\bomega}_j = \mathop{\arg\min}_{\bb_{\widehat{S}^c \backslash \{j\}}} \Big\{(2n)^{-1}\big\|\bpsi_j^{(j)} - \bpsi_{\widehat{S}^c\backslash \{j\}}^{(j)} \bb_{\widehat{S}^c\backslash \{j\}} \big\|_2^2 + \lambda_j \sum_{k \in {\widehat S}^c\backslash \{j\} } v_{jk} |b_k| \Big\}, \nonumber
\end{align}
where $\bpsi_{\widehat{S}^c \backslash \{j\}}^{(j)}$ is the matrix consisting of columns $\bpsi_{k}^{(j)}$ for $k \in \widehat{S}^c \backslash \{j\}$, $v_{jk}=\|\bpsi_k^{(j)}\|_2/ \sqrt{n}$, $b_k$ denotes the $k$th component of $\bb \in \mathbb{R}^p$ and $\lambda_j$ the regularization parameter. It is clear that $\bz_j$ is also orthogonal to the matrix $\bX_{\widehat{S} \backslash \{j\}}$ of identifiable covariates since it is a linear combination of the projected vectors $\bpsi_k^{(j)}$, $k \in \{j\} \cup \widehat{S}^c$.
Moreover, together with \eqref{zdef}, we have for any $k \in \{j\} \cup \widehat{S}^c$ that
\begin{equation}\label{zxb} 
\bz_j^\top \bpsi_k^{(j)} = \bz_j^\top (\bI_n-  \mathbb{\bP}_{\widehat{S}\backslash\{j\}})\bx_k = \bz_j^\top \bx_k.
\end{equation}
Therefore, the resulting bias factor $\eta^{\mathrm{H}}_j$ is the same as $\eta_j$ in \cite{Zhang2014}, given that
\begin{equation*}
    \eta^{\mathrm{H}}_j = \max_{k \in \widehat{S}^c \backslash \{j\}} |\bz_j^\top \bpsi_k^{(j)}|/\|\bz_j\|_2 = \max_{k \in \widehat{S}^c \backslash \{j\}} |\bz_j^\top \bx_k|/\|\bz_j\|_2 = \max_{k\ne j} |\bz_j^\top \bx_k|/\|\bz_j\|_2 = \eta_j.
\end{equation*}
In the sequel, we drop the superscript $^\mathrm{H}$ in $\eta_j ^ \mathrm{H}$ to simplify notation.
Through this two-step hybrid orthogonalization, $\bz_j$ is a strict orthogonalization against the identifiable predictors but a relaxed orthogonalization against the others.

With $(\bz_j)_{j \in [p]}$, the HOT estimator $\widehat{\bbeta} = (\hbeta_1,\ldots,\hbeta_p)^\top$ is defined coordinate-wise as
\begin{eqnarray}\label{hefei2}
\widehat{\beta}_j=\frac{\bz_j^\top\by}{\bz_j^\top\bx_j}.
\end{eqnarray}
Plugging model \eqref{eq002p} to the definition above, we obtain that
\begin{eqnarray}\label{fdrl}
\widehat{\beta}_j - \beta_j = \frac{\bz_j^\top\bveps}{\bz_j^\top\bx_j}+\sum_{k\ne j}\frac{\bz_j^\top\bx_k\beta_k}{\bz_j^\top\bx_j} = \frac{\bz_j^\top\bveps}{\bz_j^\top\bx_j}+\sum_{k \in \widehat{S}^c \backslash\{j\}}\frac{\bz_j^\top\bx_k\beta_k}{\bz_j^\top\bx_j},
\end{eqnarray}
where the last equality is due to the exact orthogonalization of $\bz_j$ to $\bX_{\widehat{S}\backslash  \{j\}}$. Furthermore, since the coefficients of order $\Omega(\sqrt{\log(p)/n})$ are typically guaranteed to be found and included in $\widehat{S}$ by the aforementioned feature screening or selection procedures under suitable conditions, the signals corresponding to the set $\widehat{S}^c\backslash \{j\}$ are weak with an aggregated magnitude no larger than the order of $s_2 \sqrt{\log (p)/n}$. Thus, by taking the relaxed projection in \emph{Step 2}, the proposed method inherits the advantage of LDPE in dealing with the weak signals. We will show that the main sparsity constraint for HOT here is that 
\[s_2 \log (p) /\sqrt{n} = o(1),\]
while $s_1$ is allowed to be  at least of order $o(n/\log p)$. 


\begin{algorithm} [htbp]
\caption{ HOT} \label{Algo} 
\scriptsize
\begin{tabbing}
   \qquad \enspace Input: $\bX \in \mathbb{R}^{n\times p}$, $\by \in \mathbb{R}^{n\times 1}$, a pre-screened set $\widehat{S}$ and a significance level $\alpha$\\
   \qquad \enspace for $j \in \{1,2,\ldots,p\}$ \\
   \qquad \qquad $\mathbb{\bP}_{\widehat{S}\setminus \{j\}} \leftarrow \bX_{\widehat{S}\setminus \{j\}}\big(\bX_{\widehat{S}\setminus \{j\}}^\top \bX_{\widehat{S} \setminus \{j\}}\big)^{-1}\bX_{\widehat{S} \setminus \{j\}}^\top$\\
   \qquad \qquad for $k \in \{j\} \cup \widehat{S}^c$ \\
   \qquad \qquad \qquad $\bpsi_k^{(j)} \leftarrow (\bI_n-\mathbb{\bP}_{\widehat{S}\backslash \{j\}})\bx_k$ \\ 
   \qquad \qquad \qquad $v_{jk}=\|\bpsi_k^{(j)}\|_2/ \sqrt{n}$ \\
   \qquad \qquad $\widehat{\bomega}_j \leftarrow \mathop{\arg\min}_{\bb_{\widehat{S}^c\backslash  \{j\}}} \Big\{(2n)^{-1}\big\|\bpsi_j^{(j)} - \bpsi_{\widehat{S}^c\backslash  \{j\}}^{(j)} \bb_{\widehat{S}^c\backslash  \{j\}} \big\|_2^2 + \lambda_j \sum_{k \in \widehat{S}^c\backslash  \{j\} } v_{jk} |b_k| \Big\}$ \\
   \qquad \qquad $\bz_j \leftarrow \bpsi_j^{(j)} - \bpsi_{\widehat{S}^c\backslash  \{j\}}^{(j)} \widehat{\bomega}_j$  \\
   \qquad \qquad $\widehat{\beta}_j \leftarrow \frac{\bz_j^\top\by}{\bz_j^\top\bx_j}$ \\
   \qquad \qquad $\tau_j  \leftarrow \|\bz_j\|_2/|\bz_j^\top \bx_j| $ \\
   \qquad \qquad $\text{CI}_j  \leftarrow [\widehat{\beta}_j - \Phi^{-1}(1-\alpha/2)\widehat{\sigma}\tau_j, \widehat{\beta}_j + \Phi^{-1}(1-\alpha/2)\widehat{\sigma}\tau_j]$ \quad \# $\widehat{\sigma}$ is the scaled lasso estimate of $\sigma$\\
   \qquad \enspace output $\{\text{CI}_j\}_{j=1}^p$ for $(1 - \alpha)$-confidence intervals of $\{\beta_j\}_{j \in [p]}$
\end{tabbing}
\end{algorithm}

The HOT inference procedure is described in Algorithm \ref{Algo}. It is interesting to notice that no initial estimate of $\bbeta$ is involved in HOT. As we shall see, HOT has two main advantages over standard de-biasing techniques. First of all, HOT completely removes the impact of identifiable signals in equation \eqref{fdrl}, thereby eliminating the bias induced by the estimation errors on identifiable coefficients and yielding more accurate confidence intervals under finite samples. Second, HOT accommodates at least $o(n/\log p)$ identifiable signals; in contrast, existing works require $s_1 = o(\sqrt{n}/\log p)$. 

\begin{figure}[htbp]
\centering
\includegraphics[height=6.3cm, width=14.5cm]{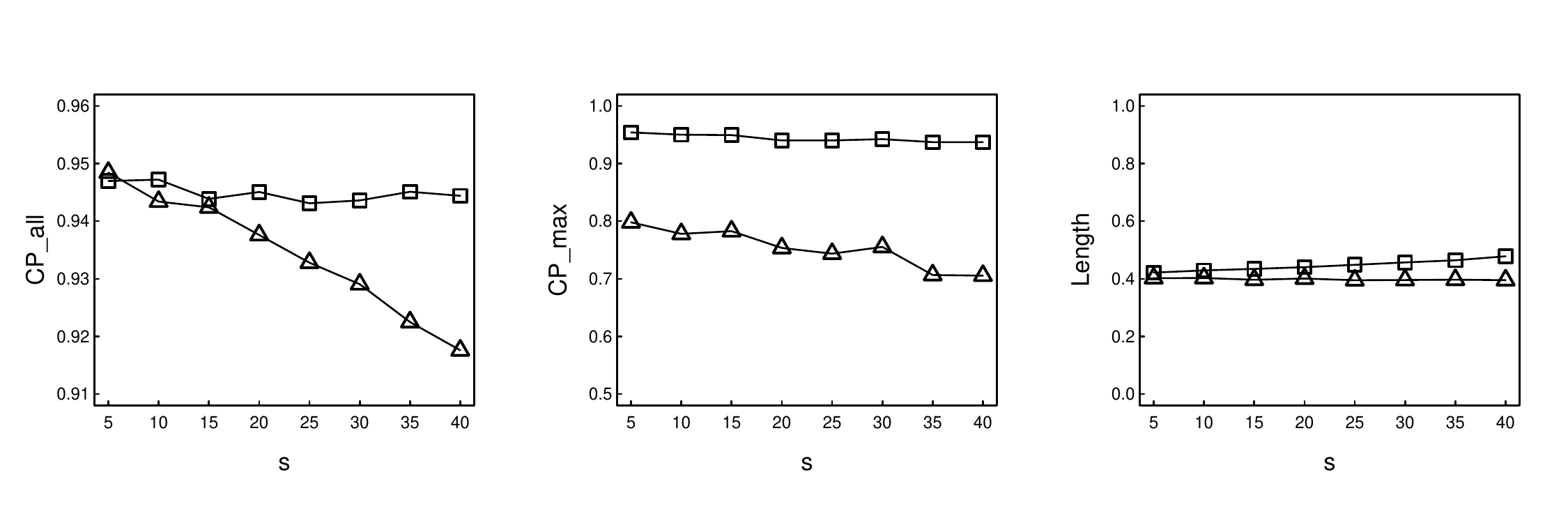}
 \caption{Inference results of HOT-SIS (square) and LDPE (triangle). CP\_all: the averaged coverage probability for all coefficients; CP\_max: the averaged coverage probability for nonzero coefficients; Length: the averaged length of confidence intervals for all coefficients; $s$: the number of nonzero coefficients.}
 \label{fg2}
\end{figure}

The enhanced efficiency of inference via HOT can be seen from a simple simulation example summarized in Figure \ref{fg2} with a pre-screened set $\widehat{S}$ selected by SIS, where the data are  generated similarly as in Section \ref{zemin} so that the number of nonzero coefficients $s$ varies from 5 to 40, the sample size $n=200$, the dimensionality $p=300$, and the covariate correlation parameter $\rho = 0.8$. It is clear that when the number of nonzero coefficients increases, the averaged coverage probability of confidence intervals by LDPE deteriorates, while the confidence intervals constructed through HOT keep the averaged coverage probability around $95\%$. 

\subsection{{An equivalent form  of the hybrid orthogonalization}}\label{edg1}
It turns out that the aforementioned two-step hybrid orthogonalization can be formulated as a partially penalized regression procedure. The following proposition elucidates this alternative formulation.

\begin{proposition}\label{pro1}
For any $j$, $1 \leq j \leq p$, the hybrid orthogonalization residual $\bz_j$ defined in (\ref{penl}) satisfies that
\begin{align*}
&\bz_j = \bx_j-\bX_{-j}\widehat{\bpi}^{(j)}_{-j},
\end{align*}
where $\widehat{\bpi}^{(j)}_{-j}$ is the minimizer of the following partially penalized regression problem:
\begin{align*}
&\widehat{\bpi}^{(j)}_{-j} = \mathop{\arg\min}_{\bb_{-j}} \Big\{(2n)^{-1}\big\| \bx_j - \bX_{-j} \bb_{-j}\big\|_2^2 +  \lambda_j \sum_{k \in \widehat{S}^c\backslash \{j\} } v_{jk} |b_k| \Big\}. \nonumber
\end{align*}
Here $\lambda_j$ and $v_{jk}$ are the same as those in problem (\ref{penl}).
\end{proposition}

Proposition \ref{pro1} shows that $\bz_j$ in Algorithm \ref{Algo} can be viewed as the residual vector from a partially penalized regression problem.
Under the context of the goodness-of-fit test for high-dimensional logistic regression, \cite{Jankova2020} exploited a similar partial penalty to better control the type I error of the test in its empirical study. \cite{Jankova2020} argued that this technique helps with controlling the reminder term of the asymptotic expansion of the test statistic under a ``beta-min condition'' that all the non-zero coefficients are at least of order $\sqrt{s \log p / n}$ in absolute value. The major difference between our work and \cite{Jankova2020} lies in the setups and goals: our motivation of imposing partial penalty is to accommodate larger sparsity of the true model in inference. Besides, an important methodological difference is that we do not residualize our predictors against the target predictor $\bx_j$ in the exact orthogonalization step if $\bx_j$ is among the identifiable predictors.

\section{Theoretical Properties}\label{tr}

This section presents the theoretical guarantees of the proposed HOT test. Suppose that $\xi$ is the smallest positive constant satisfying $\log p = O(n^\xi)$, which implies that the dimensionality $p$ can grow exponentially fast with respect to $n$. Then the maximum spurious correlation level is of order $\sqrt{\log (p)/n} = O(n^{-\frac{1 - \xi}{2}})$.
It follows that any coefficient $\beta_i$ of order $\Omega(n^{-\kappa})$ is identifiable if $\kappa < (1 - \xi)/2$. Therefore, we define $S_1$ as follows to collect all the identifiable coefficients:
$$S_1 = \{1\leq i\leq  p: |\beta_i| \geq c_1 n^{-\kappa}\},$$
where $\kappa$ can be sufficiently close to but smaller than $(1 - \xi)/2$, and where $c_1$ is some constant.
To distinguish weak signals from those in $S_1$, we assume that  $\beta_i$'s are  uniformly bounded by $O(\sqrt{(\log p)/n})$ for $i \in S_1^c$.

%

%

%
%

\subsection{Theoretical results for fixed design}
We establish the asymptotic properties of the proposed estimator under fixed design $\bX$.
\begin{theorem}\label{haolei}
Given a fixed set $ \widehat{S}\supset S_1$, the HOT estimator $\widehat{\beta}_j$ satisfies
\begin{align*}
&\widehat{\beta}_j - \beta_j=W_j+\Delta_j,\\
W_j=\frac{\bz_j^\top\bveps}{\bz_j^\top\bx_j}, \ \ |\Delta_j|\leq &\left(\max_{k \in \widehat{S}^c \backslash \{j\}}\bigg|\frac{\sqrt{n}\lambda_j \|\bpsi_k^{(j)}\|_2}{\bz_j^\top\bx_j}\bigg|\right) \sum_{k \in {{S}_1^c}}|\beta_k|.
\end{align*}
Suppose that the components of $\bveps$ are independent and identically distributed with mean 0 and variance $\sigma^2$.
If there is some $\varsigma>0$ such that the inequality
 \begin{align}\label{lya}
\lim_{n\rightarrow\infty}\frac{1}{(\sigma\|\bz_j\|_2)^{2+\varsigma}}\sum_{k=1}^n\mathbb{E}|z_{jk}\varepsilon_k|^{2+\varsigma}=0
\end{align}
holds, we further have
  \begin{align*}
\lim_{n\rightarrow\infty}\mathbb{P}(W_j\leq \tau_j\sigma x)=\Phi(x), \ \ \forall x\in \mathbb{R},
\end{align*}
where $\tau_j := \|\bz_j\|_2/|\bz_j^\top \bx_j|$, and where $\Phi(\cdot)$ is the distribution function of the standard normal distribution.

\end{theorem}

Theorem \ref{haolei} presents a preliminary theoretical result of the proposed method. The first part of this theorem shows that the error of the HOT estimator $\widehat{\beta}_j$ can be decomposed into a noise term $W_j$ and a bias term $\Delta_j$, which is consistent with the results in \cite{vandeGeer2014, Zhang2014}. However, one interesting difference from the standard de-biasing methods here is that $\Delta_j$ is non-random, because we do not need the random initial estimator $\widehat{\bbeta}^{\text{(init)}}$ as shown in \eqref{hefei2} and \eqref{fdrl}.
The proposed HOT estimator can be regarded as a modified least squares estimator for  high-dimensional settings, which is generally not shared by  these  de-biasing methods and is actually the  essential reason for the above difference.

The second part of this theorem establishes the asymptotic normality of the noise term $W_j$ based on the classical Lyapunov's condition \eqref{lya}, which  relaxes the  Gaussian assumption in \cite{vandeGeer2014, Zhang2014}.
Recall that the Gaussian assumption in \cite[Theorem 2.1]{vandeGeer2014} was imposed to guarantee the accuracy of $\|\widehat{\bbeta}^{\text{(init)}} - \bbeta\|_1$, which suggests that the initial estimator $\widehat{\bbeta}^{\text{(init)}}$ can induce extra error.

Moreover,
the assumption of independent and identically distributed noise can be further weakened to that of independent but possibly non-identically distributed noise, as detailed in the following corollary:

\begin{corollary}\label{sanxia}
Suppose that the errors $\{\varepsilon_k\}$ are independent across $k$ and that $\varepsilon_k$ has mean 0 and variance $\sigma_k^2$.
If there is some $\varsigma>0$ such that
 \begin{align}\label{lya1}
\lim_{n\rightarrow\infty}\frac{1}{A_{nj}^{2+\varsigma}}\sum_{k=1}^n\mathbb{E}|z_{jk}\varepsilon_k|^{2+\varsigma}=0
\end{align}
holds with $A_{nj}=(\sum_{k=1}^n z_{jk}^2\sigma_k^2)^{1/2}$, we further have
\begin{align*}
    \lim_{n\rightarrow\infty}\mathbb{P}\bigg(W_j\leq \frac{A_{nj}x}{\bz_j^\top\bx_j} \bigg)=\Phi(x), \ \ \forall x\in \mathbb{R}.
\end{align*}
\end{corollary}

The Lyapunov condition \eqref{lya1} here is slightly different from the version in \eqref{lya} because of the heteroskedasticity of noise.  To develop intuition for this  condition, we
consider a simple case where $\bz_j = \textbf{1}$. Then \eqref{lya1} can be rewritten as
\begin{align}\label{lya2}
    \lim_{n\rightarrow\infty}\frac{1}{(\sum_{k=1}^n\sigma_k^2)^{(2+\varsigma)/2}}\sum_{k=1}^n\mathbb{E}|\varepsilon_k|^{2+\varsigma}=0,
\end{align}
which is a pretty weak assumption about the random error $\bveps$. For example, when  $\varepsilon_k$ are independent across $k$ and distributed with the uniform distribution $U[-\sqrt{k}, \sqrt{k}]$, the above equality  \eqref{lya2} holds with $\varsigma=2$.

%
%
%
The conclusions in  Theorem \ref{haolei} are open-ended, given that the magnitude of either $|\Delta_j|$ or $\tau_j$ is not specified, and that the maximum number of identifiable signals is not explicitly derived. In the sequel, we will solidify Theorem \ref{haolei} by incorporating some random-design assumptions.
%
%

\subsection{Theoretical results for random design}
We provide here the theoretical results for random design $\bX$.
Before introducing the main results, we list a few technical conditions that underpin the asymptotic properties of the HOT test. %

\begin{condition} \label{fenbuz}
$\bX \sim N(\bzero, \bI_n \otimes \bSigma)$, where the eigenvalues of $\bPhi=\bSigma^{-1}=(\phi_{ij})_{p\times p}$ are bounded within the interval $[1/L, L]$ for some constant $L\geq 1$.
\end{condition}

\begin{condition} \label{weidu}
$s_0 (\log p)/n = o(1)$, where $s_0 = \max_{1 \leq j \leq p} s_j^*$ with $s_j^*=|\{1 \leq k \leq p: k\neq j, \phi_{jk}\neq 0\}|$.
\end{condition} 

\begin{condition} \label{xinhao} 
$s_1 =|S_1| = o(\max\{n/\log p, n/s_0\})$, and $\|\bbeta_{S_1^c}\|_1 =o(1/\sqrt{\log p})$.
\end{condition}

Condition \ref{fenbuz} assumes that the rows of the design matrix are independent and normally distributed. The eigenvalues of the precision matrix are well bounded from below and above. Similarly to \cite{Zhang2014}, the Gaussian assumption is imposed to facilitate the theoretical analysis when justifying some model identifiability condition such as the bounded sign-restricted cone invertibility factor assumption \citep{ye2010}. Our technical arguments also apply to non-Gaussian settings such as the sub-Gaussian distribution. Condition \ref{weidu} is a typical assumption in high dimensions and the sparsity of the population precision matrix is needed to guarantee consistent estimation of the penalized regression (relaxed projection) among the predictors. We will utilize this standard assumption, which was also imposed in other de-biasing works such as \cite{vandeGeer2014} and \cite{Zhang2014}, to analyze the variance of the hybrid orthogonalization  estimator.

The most appealing properties of the HOT inference are reflected in Condition \ref{xinhao}, especially its first part. To the best of our knowledge, it is weaker than the sparsity constraints in existing high-dimensional works on statistical inference or even the estimation since it holds when either $s_1 = o(n/\log p)$ or $s_1 = o(n/s_0)$. Specifically, it allows the number of identifiable nonzero coefficients $s_1 = o(\max\{n/\log p, n/s_0\})$, where $o(n/\log p)$ is a fundamental constraint on the sparsity level to guarantee consistent estimation in high dimensions. Therefore, the bound $o(n/\log p)$ is much weaker than $o(\sqrt{n}/\log p)$ imposed in the aforementioned de-biasing methods or $o(n/ (\log p)^2)$ in \cite{Javanmard2018}. Remarkably, the proposed inference method even allows for $s_1 = o(n/s_0)$, which can be substantially larger than $o(n/\log p)$ when $s_0 \ll \log p$.
The second part of Condition \ref{xinhao} implies that the sparsity level $s_2$ of the weak signals can be bounded as
\begin{equation*}
s_2 = \sum_{j \in S_1^c} \min{\big\{|\beta_{j}|/(\sigma \lambda_0), 1 \big\}} \asymp \|\bbeta_{S_1^c}\|_1/(\sigma \lambda_0) = o\bigg(\frac{1}{\sqrt{\log p}}  \frac{\sqrt{n}}{\sqrt{\log p}}\bigg) = o(\sqrt{n}/\log p),
\end{equation*}
where we recall $\lambda_0 = \sqrt{2 \log p / n}$.
It is worth mentioning that when $s_1 \gg n/\log p$, the identifiable coefficients may not be found by feature selection methods due to the lack of estimation consistency. Feature screening procedures such as SIS, distance correlation learning, and HOLP will be preferred. 

We first  analyze the statistical properties of the relaxed projection in (\ref{penl}), which is essentially a lasso regression on the residualized features $\bpsi_{\widehat{S}^c \backslash \{j\}}^{(j)}$. We start with investigating the sign-restricted cone invertibility factor \citep{ye2010} of $\bpsi_{\widehat{S}^c \backslash \{j\}}^{(j)}$, an important quantity that determines the error of the lasso estimator.
Let $S_{j*}=\{1 \leq k \leq p: k\in {\widehat{S}^c\backslash \{j\}}, \phi_{jk}\neq 0\}$ and $S_{j*}^c=\widehat{S}^c\backslash (\{j\}\cup S_{j*})$. For some constant $\xi\geq 0$, the sign-restricted cone invertibility factor of the design matrix $\bpsi_{\widehat{S}^c\backslash \{j\}}^{(j)}$ is defined as
\begin{align*}
  F_2^{(j)} (\xi,S_{j*}) := \inf \left\{|S_{j*}|^{1/2}\Big\|n ^ {-1}(\bpsi_{\widehat{S}^c\backslash \{j\}}^{(j)})^\top \bpsi_{\widehat{S}^c\backslash\{j\}}^{(j)} \bu\Big\|_{\infty}/ \|\bu\|_{2}: \bu \in \mathcal{G}_{-}^{(j)}(\xi,S_{j*}) \right\}, 
\end{align*}
where the sign-restricted cone 
$$\mathcal{G}_{-}^{(j)}(\xi,S_{j*})=\Big\{\bu:\|\bu_{S_{j*}^c}\|_1 \leq \xi \|\bu_{S_{j*}}\|_1, \ u_k(\bpsi_k^{(j)})^\top\bpsi_{\widehat{S}^c\backslash \{j\}}^{(j)} \bu/n\leq 0, \ \forall \ k\in S_{j*}^c\Big\}.$$
Compared with the compatibility factor or the restricted eigenvalue \citep{Bickel09}, the sign-restricted cone invertibility factor can be more accommodating. In fact, all three of them are satisfied with high probability when the population covariance of the random design has bounded spectrum. The following proposition shows that the sign-restricted cone invertibility factor of the design matrix $\bpsi_{\widehat{S}^c\backslash\{j\}}^{(j)}$ is well above zero with high probability. 
It ensures the identifiability of the support of the true regression coefficients by avoiding collinearity among the design vectors when constrained on the sign-restricted cone.

\begin{proposition}\label{pro2}
When Conditions \ref{fenbuz} and \ref{weidu} hold and $d = o(n/s_0)$, we have for any $j$, $1 \leq j \leq p$, and sufficiently large $n$, with probability at least $1- o(p^{1-\delta})$, there exists a positive constant $c^*$ such that
\begin{align*}
 F_2^{(j)} (\xi,S_{j*})  \geq c^*.
\end{align*}
\end{proposition}

Note that the conclusion of Proposition \ref{pro2} is typically assumed directly in high dimensions for design matrices with independent and identically distributed (i.i.d.) rows. However, the design matrix $\bpsi_{\widehat{S}^c\backslash\{j\}}^{(j)}$ here consists of projected vectors, so that the rows of $\bpsi_{\widehat{S}^c\backslash \{j\}}^{(j)} = (\bI_n-\mathbb{\bP}_{\widehat{S}\backslash \{j\}})\bX_{\widehat{S}^c\backslash \{j\}}$ are neither independent nor identically distributed.
We overcome this difficulties by a new technical argument that is tailored for the $L_1$-constrained structure of the cone. The upper bound $d = o(n/s_0)$ is required to show the asymptotically vanishing difference between $\bpsi_{\widehat{S}^c\backslash \{j\}}^{(j)}$ and its population projected counterpart.

Proposition \ref{pro2} lays down the foundation for the following theorem on the statistical error of the relaxed projection.
\begin{theorem}\label{theob}
Let $\lambda_j=(1+\epsilon)\sqrt{2\delta \log(p)/(n\phi_{jj})}$ for some constants $\epsilon>0$ and $\delta > 1$ and any $1\leq j \leq p$. Given a set $\widehat{S}$, when Conditions \ref{fenbuz} and \ref{weidu} hold and $d = |\widehat{S}| = o(n/s_0)$, we have for any $1 \le j \le p$ and sufficiently large $n$ that with probability at least $1- o(p^{1-\delta})$,
\begin{align*}
& \|\widehat{\bomega}_j-{\bomega}_j\|_q = O(s_0^{1/q}\lambda_j) = O\{s_0^{1/q} \sqrt{\log(p)/n}\}, \quad q \in [1,2];\\
& n^{-1/2} \|\bpsi_{\widehat{S}^c \backslash \{j\}}^{(j)}(\widehat{\bomega}_j-{\bomega}_j)\|_2 = O(s_0^{1/2}\lambda_j) = O\{\sqrt{s_0\log(p)/n}\},
\end{align*}
where ${\bomega}_j=-\bPhi_{j,{\widehat{S}^c \backslash \{j\}}}^\top/\phi_{jj}$ is the population counterpart of $\widehat{\bomega}_j$. 
\end{theorem}

Theorem \ref{theob} establishes estimation and prediction error bounds for the relaxed projection in (\ref{penl}). One striking feature here is that the statistical rates no longer depend on the size of the penalization-free set $\widehat{S}$: the identifiable predictors have been eliminated in the exact orthogonalization step and are not involved in the relaxed projection.
Next, we analyze the partially penalized regression, which is crucial to derive the asymptotic normality of the subsequent HOT test.

\begin{theorem}\label{theoa}
Let $\lambda_j=(1+\epsilon)\sqrt{2\delta \log(p)/(n\phi_{jj})}$ for some constants $\epsilon>0$ and $\delta > 1$, $1 \leq j \leq p$. Given any  pre-screened set $\widehat{S}$, when Conditions \ref{fenbuz} and \ref{weidu} hold and $d = |\widehat{S}| = o\{n/\log (p)\}$, for any $1 \leq j \leq p$ and sufficiently large $n$, we have with probability at least $1- o(p^{1-\delta})$ that
\begin{align*}
& \|\widehat{\bpi}_{-j}^{(j)}-\bpi_{-j}^{(j)}\|_q = O\{(s_0+d)^{1/q}\lambda_j\} = O\{(s_0+d)^{1/q}\sqrt{\log(p)/n}\}, \forall q \in [1,2];\\
& n^{-1/2} \|\bX_{-j} (\widehat{\bpi}_{-j}^{(j)}-\bpi_{-j}^{(j)})\|_2 = O\{(s_0+d)^{1/2}\lambda_j\} = O\{\sqrt{(s_0+d)\log(p)/n}\},
\end{align*}
where $\bpi_{-j}^{(j)}=-\bPhi_{j,-j}^\top/\phi_{jj}$ is the population counterpart of $\widehat{\bpi}_{-j}^{(j)}$. 
\end{theorem}

Theorem \ref{theoa} establishes estimation and prediction error bounds for the solution of the partially penalized regression problem. It guarantees estimation and prediction consistency of $\widehat\bpi^{(j)}_{-j}$ when the size of the penalization-free set $\widehat{S}$ is of order $o\{n/\log (p)\}$. Such an $\widehat S$ can typically be obtained through regularization methods, e.g, \citep{Bickel09, FanLv2014, Kong16}.
The estimation and prediction consistency of the solution of the partially penalized regression is important for ensuring valid inference via HOT in that it guarantees the convergence of the hybrid orthogonalization vector $\bz_j$ to its population counterpart $\textbf{e}_j = \bx_j - \mathbb{E}(\bx_j|\bX_{-j})$.
Based on Proposition \ref{pro1} and a joint analysis of penalized and penalization-free covariates, Theorem \ref{theoa} also explicitly characterizes the dependence of the estimation error on the size of $\widehat{S}$.

It is worth pointing out that the partially penalized procedure here is essentially different from those weighted penalization methods such as the adaptive lasso, whose weights target at improving the estimation accuracy of the regression coefficients. In contrast, the partial penalization in our setup can downgrade the estimation quality. Note that the pre-screened set $\widehat{S}$ of identifiable predictors may have little overlap with the support of $\bpi_{-j}^{(j)}$, which is the set of significant covariates associated with the predictor $\bx_j$. Therefore, lifting regularization on $\widehat S$ incurs more exposure of $\widehat \bpi_{-j} ^ {(j)}$ to noise and increases statistical error: the sparsity factor in the error bounds in Theorem \ref{theoa} is $s_0 + d$ instead of just $s_0$. Despite of the sacrifice on the estimation accuracy, the resulting $\bz_j$ is strictly orthogonal to $\bX_{\widehat S}$, which is the key to accommodating more identifiable signals in the HOT test than in the standard debiasing-based test.

 {In fact, both Theorems   \ref{theob} and \ref{theoa} lead to the same guarantee for $\bz_j$:}
\begin{eqnarray*}
{\mathbb{P}\left\{\|\bz_j\|_2\asymp n^{1/2}\right\}\geq 1-  o(p^{1-\delta})};
\end{eqnarray*}
{see the proof of Theorem  \ref{theo1} for details. Perhaps more importantly, the two theorems complement each other in terms of the requirement on the magnitude of $|\widehat{S}|$, so that we allow $|\widehat{S}| =  o\{\max(n/\log p, o(n/s_0))\}$.
Therefore, if $\log (p) \gg s_0$, we follow Theorem \ref{theob} to analyze $\bz_j$; otherwise we follow Theorem \ref{theoa}.}


Before presenting the main theorem on the HOT test, the following definition characterizes the properties that a pre-screened set $\widehat{S}$ should satisfy. 
\begin{definition}[Acceptable set]\label{C2}
Set $\widehat{S}$ is called an acceptable pre-screened set if it satisfies: (1) $\widehat{S}$ is independent of the data $(\bX, \by)$; 
(2) with probability at least $1-\theta_{n,p}$ for some asymptotically vanishing $\theta_{n,p}$, $S_1 \subset\widehat{S}$ and $d = |\widehat{S}| = O(s_1)$. 
\end{definition}
The first property in Definition \ref{C2} is imposed so that the randomness of $\widehat{S}$ does not interfere with the subsequent inference procedure. The same property was suggested in \cite{Fan13,Fan15} to facilitate the theoretical analysis. In practice, it can be ensured through data splitting and our numerical studies show that it can be more of a technical assumption than a practical requisite. The second property is the sure screening property for the identifiable predictors, which can be satisfied by the aforementioned screening or variable selection procedures under suitable conditions. For example, under the assumptions of bounded $\|\bbeta\|_2$ and significant marginal correlations, SIS \citep[Theorem 3]{Fan2008} can yield a sure screening set with failure probability $\theta_{n,p}$ of order $O(\exp\{-n^{1 - 2\kappa}/\log n\})$. The tail probability also vanishes asymptotically for variable selection methods \citep{Zhang10} when the design matrix satisfies the irrepresentable condition \citep{zhao06}. 
In addition, the second part of Definition \ref{C2} also assumes the size of the pre-screened set $\widehat{S}$ to be of the same order as $s_1$, the number of the true identifiable predictors. For feature screening procedures such as SIS, the selected model can be reduced to a desirable size through multiple iterations. For variable selection methods such as the lasso, the number of selected features is generally around $O(s_1)$ with probability tending to one \citep{Bickel09} if the largest eigenvalue of the Gram matrix is bounded from above. Other nonconvex regularization methods tend to select even sparser models \citep{FanLv2014,Kong16}. In fact, even if this assumption is violated, our theoretical results still hold as long as $d = o(\max\{n/\log p, n/s_0\})$ in view of the technical arguments.


Now we are ready to present our main theorem that establishes the asymptotic normality of the HOT test statistic.

\begin{theorem}\label{theo1}
Let $\lambda_j=(1+\epsilon)\sqrt{2\delta \log(p)/(n\phi_{jj})}$ for some constants $\epsilon>0$ and $\delta > 1$, $1\leq j \leq p$. Assume that Conditions \ref{fenbuz}-\ref{xinhao} hold and that $\widehat{S}$ satisfies Definition \ref{C2}. Then for any $1 \leq j \leq p$ and sufficiently large $n$, the HOT estimator $\widehat{\beta_j}$ satisfies
\begin{align*}
    \widehat{\beta}_j - \beta_j=&W_j+\Delta_j,
\end{align*}
where $W_j=\frac{\bz_j^\top\bveps}{\bz_j^\top\bx_j}$ and $|\Delta_j|=o_{\mathbb{P}}(n^{-1/2})$.
Moreover, it holds with probability at least $1-o(p^{1-\delta})-\theta_{n,p}$ that
$$\lim_{n \to \infty} \tau_j n^{1/2} = \phi_{jj}^{-1/2}.$$
Suppose that the components of $\bveps$ are independent and identically distributed with mean 0 and variance $\sigma^2$.
For any given $\bX$, if there is some $\varsigma>0$ such that the inequality
 \begin{align*}
\lim_{n\rightarrow\infty}\frac{1}{(\sigma\|\bz_j\|_2)^{2+\varsigma}}\sum_{k=1}^n\mathbb{E}|z_{jk}\varepsilon_k|^{2+\varsigma}=0
\end{align*}
holds, we further have
  \begin{align*}
\lim_{n\rightarrow\infty}\mathbb{P}(W_j\leq \tau_j\sigma x)=\Phi(x), \ \ \forall x\in \mathbb{R},
\end{align*}
 where $\Phi(\cdot)$ is the distribution function of the standard normal distribution.
\end{theorem}

Compared with existing inference results on de-biasing, Theorem \ref{theo1} guarantees the asymptotic normality of the HOT estimator under a much relaxed constraint on the number of strong signals.
When $s_0 \gg \log p$ so that $s_1 = o(n/\log p)$, either feature screening or selection methods can be used to select $\widehat S$. When $s_0 \ll \log p$ and thus $s_1 = o(n/s_0)$, feature screening procedures will be preferred to identify $\widehat S$. The constraint $s_1 s_0 = o(n)$ in the latter case reveals an interesting trade-off between the sparsity in the regression coefficient and that in the precision matrix of the predictors. A recent related work \cite{Javanmard2018} used a novel leave-one-out argument to show that when $s_0 = o(\sqrt{n}/ \log p)$, the standard de-biasing test allows $s = o(n/ (\log p)^2)$, which is still more stringent than $s_1 = o(\max(n / \log p, n / s_0))$ in our result. Finally, while coping with much more strong signals, the HOT test also accommodates $o(\sqrt{n}/\log p)$ weak signals as the standard de-biasing tests do.  

Moreover, $n ^ {1 /2}\tau_j$ converges to the same constant as that of LDPE in \cite{Zhang2014}. Given that the noise factor determines the asymptotic variance of the test, Theorem \ref{theo1} implies no efficiency loss incurred by HOT.
Similarly to \cite{Javanmard2014} and \cite{vandeGeer2014}, when the noise standard deviation $\sigma$ is unknown in practice, we can replace it with some consistent estimate $\widehat{\sigma}$ given by the scaled lasso for general Gaussian settings, say. When the data are heteroscedastic, we can use the kernel method  as suggested in \cite{muller1987estimation}  to estimate $\sigma$.

\section{Simulation Studies}\label{sr}
In this section, we use simulated data to investigate the finite sample performance of HOT, in comparison with that of LDPE. Both SIS and HOLP are utilized to pre-screen the set $\widehat{S}$ of identifiable coefficients and we denote the two versions of HOT by HOT-SIS and HOT-HOLP, respectively. To control the number of selected identifiable coefficients in $\widehat{S}$, we take the least squares estimates after constraining the covariates on $\widehat{S}$ and then apply BIC to choose the optimal one. While these two methods use the same data set for pre-screening and the subsequent inference procedure, as a comparison, an independent data set of the same size as the original sample is also generated for the pre-screening step. The corresponding versions are denoted by HOT-SIS(I) and HOT-HOLP(I), respectively. 

The hybrid orthogonalization vectors $\bz_j$ are calculated by the lasso as defined in \eqref{penl} with the regularization parameter tuned by GIC \citep{Fan2013}.
Two simulation examples are conducted for both sparse and non-sparse settings. Throughout the simulation studies, we aim at establishing the $95\%$ confidence intervals for all coefficients, corresponding to the significance level $\alpha = 0.05$. The estimated error standard deviation $\widehat{\sigma}$ is obtained through the scaled lasso with the universal regularization parameter set according to \cite{Sun13}. 



\subsection{Simulation Example 1}\label{zemin}

Here we adopt the same setting as that in \cite{vandeGeer2014}. To be specific, 100 data sets were simulated from model (\ref{eq002p}) with $(n,p,\sigma)=(100,500,1)$ so that $\bveps\sim N(\bzero, \bI_n)$. For each data set, the rows of $\bX$ were sampled as independent and identically distributed copies from $N(\bzero,\bSigma)$, where $\bSigma=(\rho^{|j-k|})_{p\times{p}}$ with $\rho = 0.9$. The true coefficient vector $\bbeta = (\beta_1, \beta_2,..., \beta_{s_0}, 0,...,0)^\top$, where $s_0=15$ and the nonzero $\beta_j$ were independently sampled from the uniform distribution $U[0, 2]$.
Three performance measures are used to evaluate the inference results, including the averaged coverage probability for all coefficients (CP\_all), the averaged coverage probability for nonzero coefficients (CP\_max) in this example or that for the maximum coefficients in the second example, and the averaged length of confidence intervals for all coefficients (Length). The averaged estimated error standard deviation $\widehat{\sigma}$ obtained by the scaled lasso is $0.95$ and we also display the performance measures when $\widehat{\sigma}$ equals to the population value $1$ for comparison. 

%
%
%
%

\begin{table}[H]
\centering
\caption{\label{my-label2} Three performance measures by different methods over 100 replications in Section \ref{zemin}.}
\medskip
\begin{tabular}{ccccccc}
\hline
					&         &HOT-SIS &HOT-SIS(I)&HOT-HOLP &HOT-HOLP(I)&LDPE       \\
\hline
&CP\_all                                   & 0.942        & 0.943       & 0.941       & 0.939    & 0.943        \\
$\widehat{\sigma} = 0.95$ &CP\_max         & 0.942        & 0.945       & 0.935       & 0.941    & 0.702      \\
&Length      	                           & 0.673        & 0.677       & 0.670       & 0.671    & 0.659       \\

&CP\_all                                   & 0.951        & 0.953       & 0.950       & 0.951    & 0.953        \\
$\widehat{\sigma} = 1$ &CP\_max            & 0.950        & 0.951       & 0.947        & 0.951    & 0.729      \\
&Length      	                           & 0.700        & 0.701       & 0.701       & 0.700    & 0.697      \\
\hline
\end{tabular}
\end{table}

%
%
%
%
The results are summarized in Table \ref{my-label2}. It is clear that while the averaged coverage probabilities for all coefficients (CP\_all) by different methods are all around $95\%$, the averaged coverage probabilities for nonzero coefficients (CP\_max) of different versions of HOT are significantly closer to the target than that of LDPE under this finite sample set. Moreover, the averaged lengths of confidence intervals by different methods are around the same level.
By comparing the results of HOT using the same or an independent data set for the pre-screening step, we find that the independence assumption on the pre-screening set can be more of a technical assumption than a practical necessity. 

\subsection{Simulation Example 2}\label{sim2}
In this second example, we consider the approximately sparse setting in \cite{Zhang2014}, where all coefficients are nonzero but decay fast as the indices increase except for some significant ones. Specifically, the true regression coefficients $\beta_j = 3\lambda_{\rm{univ}}$ with $\lambda_{\rm{univ}} = \sqrt{2(\log p)/n}$ for $j= 200, 400, \dots, 1000$, and $\beta_j = 3\lambda_{\rm{univ}}/j^2$ for all the other $j$. The other setups are similar to those in the first example with $(n,p,\sigma, \rho)=(200,1000,1,0.5)$. We adopt the same performance measures and summarize the results in Table \ref{my-label1}. Under such non-sparse setting, the scaled lasso estimate $\widehat{\sigma}$ tended to overestimate $\sigma$ a bit with its mean value equaling $1.12$. Therefore, with the estimated $\widehat{\sigma}$, the averaged coverage probabilities of different versions of HOT are above the target $95\%$ a bit, while that of LDPE for maximum coefficients is still slightly below $95\%$. If we utilize the population value so that $\widehat{\sigma} = 1$ to construct the confidence intervals, the averaged coverage probabilities of different versions of HOT are around the target. It shows that the proposed method can still be valid under some approximately sparse settings. 

\begin{table}[H]
\centering
\caption{\label{my-label1}Three performance measures by different methods over 100 replications in Section \ref{sim2}.}
\medskip

\begin{tabular}{ccccccc}
\hline
					&         &HOT-SIS &HOT-SIS(I)&HOT-HOLP &HOT-HOLP(I)&LDPE       \\
\hline
&CP\_all                                   & 0.971        & 0.970       & 0.970       & 0.967    & 0.949        \\
$\widehat{\sigma} = 1.12$ &CP\_max         & 0.967        & 0.967       & 0.970       & 0.965    & 0.935      \\
&Length      	                           & 0.341        & 0.342       & 0.340       & 0.338    & 0.341       \\
&CP\_all                                   & 0.949        & 0.950       & 0.949       & 0.947    & 0.922        \\
$\widehat{\sigma} = 1$ &CP\_max            & 0.943        & 0.945       & 0.945       & 0.943    & 0.896      \\
&Length      	                           & 0.305        & 0.305       & 0.305       & 0.305    & 0.306      \\
\hline
\end{tabular}
\end{table}

%
%
%
\section{Application to Stock Short Interest and Aggregate Return Data}\label{rd} 
In this section, we will demonstrate the usefulness of the proposed methods by analyzing the stock short interest and aggregate return data set originally studied in \cite{Rapach16}, available at Compustat (www.compustat.com). The raw data set includes 3274 firms' monthly short interests, reported as the numbers of shares that were held short in different firms, as well as the $\text{S}\&\text{P}$ 500 monthly log excess returns from January 1994 to December 2013. As shown in \cite{Rapach16}, the short interest index, that is, the average of short interests of the firms, is arguably the strongest predictor of aggregate stock returns. Then with the proposed large-scale inference technique, we can further identify which firms among the thousands of candidates that had significant impacts on the $\text{S}\&\text{P}$ 500 log excess returns.

We approach this problem by the following predictive regression model 
 \begin{equation}\label{ML1}
   r_{t+1}= \bx_t^\top \bbeta +\epsilon_{t+1},
 \end{equation}
where $r_{t + 1}$ is the $\text{S}\&\text{P}$ 500 log excess return of the $(t + 1)$th month and $\bx_t=(x_{1t},\ldots,x_{pt})^\top$ with each component $x_{jt}$ indicating the short interest of the $j$th firm in the $t$th month.
Due to the high correlations between short interests of firms, we use one representative for those firms with correlation coefficients larger than 0.9, resulting in $p = 975$ representative firms. Moreover, the aforementioned period yields a sample size of $n=240$ and both predictor and response vectors are standardized to have mean zero and $L_2$-norm $\sqrt{n}$. We applied HOT and LDPE similarly as in Section \ref{sr} at a significance level of $\alpha=0.005$. SIS was utilized in HOT for screening the identifiable predictors to bypass the high correlations.
\begin{table}[H]
\centering
\caption{\label{my-label33} Numbers of important firms and averaged lengths of confidence intervals at a significance level of $\alpha=0.005$ by different methods in Section \ref{rd}.}
\medskip
\begin{tabular}{ccccc}
\hline
                  &                      & HOT         & LDPE    \\ 
\hline
           &Num                             & 6          & 11        \\ 
          &Length                           & 0.489     & 0.507    \\  
          \hline
\end{tabular}
\end{table}
Table \ref{my-label33} summarizes numbers of important firms and averaged lengths of confidence intervals by different methods.
While the averaged lengths of confidence intervals are very close, LDPE selects more firms than HOT. To verify the significance of the selected firms, we refit model (\ref{ML1}) via ordinary least squares (OLS) constrained on the selected firms by HOT and LDPE, respectively. The corresponding significance test shows that only five firms out of the eleven selected by LDPE are significant at the 0.05 level, while four our of the six firms identified by HOT are significant around the 0.05 level. We display the refitted OLS estimates and the p-values for the six firms by HOT in Table \ref{mytable333}. 

\begin{table}[H]
  \centering
   \caption{Refitted OLS estimates $\widehat{\beta}_j$ and p-values for the six firms selected by HOT.}\label{mytable333}
   \smallskip
  \begin{tabular}{lcc|lcc}
\hline
     Firm & $ \widehat{\beta}_j$ & $p$-value & Firm & $\widehat{\beta}_j$ & $p$-value \\
     \hline
     GAPTQ  & -0.167       & 0.025     & PKD   & -0.098      &0.183  \\
     WHR    & 0.088        & 0.361     & HAR   & 0.171       &0.034  \\
     GCO    & -0.310       & 0.000     & OMS   & -0.126      &0.052  \\
     \hline
     \end{tabular}
\begin{threeparttable}
\begin{tablenotes}
        \footnotesize
        \item[1] GAPTQ: Great Atlantic $\&$ Pacific Tea Company; WHR: Whirlpool Corporation; GCO: Genesco Inc; PKD: Parking Drilling Company; HAR: Harman International Industries Inc;
        OMS: Oppenheimer Multi-Sector
      \end{tablenotes}
    \end{threeparttable}
\end{table}

In view of Table \ref{mytable333}, the p-values of four firms, GAPTQ, GCO, HAR, and OMS, are around or below 0.05. These firms are also influential in their respective fields. Among them, GAPTQ (Great Atlantic $\&$ Pacific Tea Company) is a supermarket with estimated sales of $\$6.7$ billion in 2012, ranking 28th among the ``top 75 food retailers and wholesalers''. The most significant firm GCO (Genesco Inc) is a retailer of branded footwear and clothing with thousands of retail outlets across the world. Furthermore, OMS (Oppenheimer Multi-Sector) has grown into one of the largest and the most popular investment managers in the world since going public in 1959. The negative coefficients of them imply that the increases in the short interests of these firms tend to significantly decrease the future $\text{S}\&\text{P}$ 500 return. With these identified influential firms, it can be helpful to develop some augmented short interest index to further enhance the prediction performance for excess returns. 

\section{Discussions}\label{dis} 
In this paper, we have demonstrated that the maximum number of identifiable coefficients is allowed to be substantially larger than existing ones in presence of many possibly nonzero but weak signals for valid high-dimensional statistical inference. We also propose the hybrid orthogonalization  estimator that takes full advantage of the identifiable coefficients, so that the bias induced by the estimation errors of the identifiable coefficients can be removed to enhance the inference efficiency. 


Alternatively, based on the hybrid orthogonalization vectors $\bz_j$ ($1 \leq j \leq p$) constructed in Section \ref{ms}, we can also utilize the self-bias correction idea with the aid of an initial estimate such as the scaled lasso estimate $\widehat{\bbeta}^{\text{(init)}}=(\hbeta_1^{\text{(init)}},\ldots,\hbeta_p^{\text{(init)}})^\top$ in equation \eqref{sl}. Then an alternative hybrid orthogonalization  estimator $\widehat{\bbeta}^{\rm A} = (\hbeta_1^{\rm A}, \cdots, \hbeta_p^{\rm A})^\top$ can be defined through each coordinate as 
\begin{equation*}
\hbeta_j^{\rm A} = \hbeta_j^{\text{(init)}} + \bz_j^\top (\by - \bX \widehat{\bbeta}^{\text{(init)}})/\bz_j^\top \bx_j.
\end{equation*}

The validity of it can be seen from the following decomposition
\begin{equation*}
\hbeta_j^{\rm A} - \beta_j = \frac{\bz_j^\top\bveps}{\bz_j^\top\bx_j} + \sum_{k\ne j}\frac{\bz_j^\top\bx_k(\beta_k-\widehat{\beta}_k^{(\rm init)})}{\bz_j^\top\bx_j} = \frac{\bz_j^\top\bveps}{\bz_j^\top\bx_j} + \sum_{k \in \widehat{S}^c\backslash\{j\}} \frac{\bz_j^\top\bx_k(\beta_k-\widehat{\beta}_k^{(\rm init)})}{\bz_j^\top\bx_j},
\end{equation*}
where the items with indices in the set $\widehat{S}\backslash  \{j\}$ are zero due to the exact orthogonalization of $\bz_j$ to $\bX_{\widehat{S}\backslash  \{j\}}$ similarly as in \eqref{fdrl}. Moreover, the bias term on the right hand side of the above equality is bounded by
\begin{align*}
\sum_{k \in \widehat{S}^c\backslash\{j\}}\frac{|\bz_j^\top\bx_k(\beta_k-\widehat{\beta}_k^{(\rm init)})|}{|\bz_j^\top\bx_j|} &\leq \big(\max_{k\ne j} \frac{|\bz_j^\top \bx_k|}{\|\bz_j\|_2}\big) \cdot \frac{\|\bz_j\|_2}{|\bz_j^\top \bx_j|} \cdot \|(\widehat{\bbeta}^{\text{(init)}} - \bbeta)_{\widehat{S}^c\backslash\{j\}}\|_1 \\
&= \eta_j \tau_j \|(\widehat{\bbeta}^{\text{(init)}} - \bbeta)_{\widehat{S}^c\backslash\{j\}}\|_1, 
\end{align*}
where the bias factor $\eta_j$ and the noise factor $\tau_j$ are of the same definitions as those in \eqref{sbc} except that the orthogonalization vectors $\bz_j$ are now the hybrid ones.

Since $\eta_j$ and $\tau_j$ can be shown to take the same orders of magnitudes as those in \cite{Zhang2014}, the main advantage of the alternative hybrid orthogonalization  estimator $\widehat{\bbeta}^{\rm A}$ is to reduce the item $\|\widehat{\bbeta}^{\text{(init)}} - \bbeta\|_1$ in the bias term of \eqref{sbc} to $\|(\widehat{\bbeta}^{\text{(init)}} - \bbeta)_{\widehat{S}^c\backslash\{j\}}\|_1$, which denotes the $L_1$-estimation error of $\bbeta$ constrained on the index set $\widehat{S}^c\backslash\{j\}$. When the possibly nonzero but weak signals are not identified by the initial estimate $\widehat{\bbeta}^{\text{(init)}}$, the bias term here would be as small as that for the proposed hybrid orthogonalization  estimator $\widehat{\bbeta}$, so that their sample size requirements for valid inference are basically the same. 


Note that except for the linear regression models that utilize the squared error loss, the second order derivative of the loss function would typically rely on the regression coefficient vector. Therefore, the main advantage of the alternative hybrid orthogonalization  estimator $\widehat{\bbeta}^{\rm A}$ is that it can be readily extended to other convex loss functions such as generalized linear models. The key point is to construct hybrid orthogonalization  vectors based on the second order derivatives of the loss function for the corresponding models after identifying the significant coefficients by powerful feature screening techniques \citep{Li2012, Ma17, Liu20} targeted at non-linear models. It would be interesting to explore the theoretical properties of hybrid orthogonalization  estimators in these extensions for future research.

\section*{Acknowledgements}
%
%
%

This work was supported by National Natural Science Foundation of China (Grants 72071187, 11671374, 71731010, 12101584 and 71921001), Fundamental Research Funds for the Central Universities (Grants WK3470000017 and WK2040000027), and China Postdoctoral Science Foundation (Grants 2021TQ0326 and 2021M703100). Yang Li and Zemin Zheng are co-first authors. 

\newpage
\setcounter{page}{1}
\setcounter{section}{0}
\setcounter{equation}{0}

\renewcommand{\theequation}{A.\arabic{equation}}
\setcounter{equation}{0}
\bigskip
\begin{center}
{\large\bf E-companion  to ``High-dimensional inference via hybrid orthogonalization''}
\end{center}

\smallskip
\begin{center}

Yang Li, Zemin Zheng, Jia Zhou and Ziwei Zhu   

\end{center}

\medskip

\noindent This supplementary material consists of two parts. Section \ref{EC1} lists the key lemmas and presents the
proofs for main results. Additional technical proofs for the lemmas are provided in Section \ref{EC2}.


\section{Proofs of main results}\label{EC1}
\subsection{Lemmas}
The following lemmas will be used in the proofs of the main results.
\begin{lemma}\label{tongfb}
Suppose that $\bX \sim N(0, \bI_n \otimes \bSigma)$ and the eigenvalues of $\bPhi=\bSigma^{-1}=(\phi_{ij})_{p\times p}$ are bounded within the interval $[1/L, L]$ for some constant $L\geq 1$. For a given set $\widehat{S}$ and any $j\in \widehat{S}^c$, we have the following conditional distribution
\begin{align*}
& \bx_j=\bX_{\widehat{S}} \bgamma_j^{(j)}+ \brho_j^{(j)},
\end{align*}
where the residual vector $\brho_j^{(j)}=\bx_j-\mathbb{E}(\bx_j|\bX_{\widehat{S}})$ and $\bgamma_j^{(j)}$ is the corresponding regression coefficient vector. 
Then the residual matrix $\brho^{(j)}_{\widehat{S}^c\backslash \{j\}}$ consisting of the residual vectors $\brho_j^{(j)}$ for $j \in \widehat{S}^c\backslash \{j\}$ satisfies
$$\brho^{(j)}_{\widehat{S}^c\backslash \{j\}}\sim N(\bzero, \bI_n \otimes \bPhi_{\widehat{S}^c\backslash \{j\}\widehat{S}^c\backslash \{j\}}^{-1}).$$
\end{lemma}
\begin{lemma}\label{sigma}
Suppose that the eigenvalues of a symmetric matrix $\bPhi=(\phi_{ij})_{p\times p}=\bSigma^{-1}=(\sigma_{ij})^{-1}_{p\times p}$ are bounded within the interval $[1/L, L]$ for some constant $L\geq 1$. Then we have
\begin{enumerate}
\item the eigenvalues of $\bSigma$ are also bounded within the interval $[1/L, L]$;
\item the diagonal components of $\bPhi$ and $\bSigma$ are bounded within the interval $[1/L, L]$;
\item the eigenvalues of any principal submatrix of $\bPhi$ and $\bSigma$ are bounded within the interval $[1/L, L]$;
\item the $L_2$-norms of the rows of $\bPhi$ and $\bSigma$ are bounded within the interval $[1/L, L]$;
\item  $\sigma_{jj}\phi_{jj}\geq 1$ for any $j\in \{1,\ldots,p\}$.
\end{enumerate}

\end{lemma}


\begin{lemma}\label{lse-lem}
Consider the following linear regression model
\begin{eqnarray*}\label{lse0}
\by=\bX\bbeta+\bveps, \quad  \bveps\sim N(\bzero,\sigma^2 \bI_n),
\end{eqnarray*}
where the $n\times d$ design matrix  $\bX$ is independent of $\bveps$. Suppose that $\bX \sim N(\bzero, \bI_n \otimes \bSigma)$ and the eigenvalues of $\bSigma=(\sigma_{ij})_{d\times d}$ are bounded within the interval [1/L, L] for some constant $L\geq 1$. When $d = o(n)$, for any $p$ satisfying $\log(p)/n=o(1)$ and any constant $\delta > 1$, we have with probability at least $1- o(p^{1-\delta})$, the ordinary least squares estimator $\widehat{\bbeta}=(\bX^\top\bX)^{-1}\bX^\top\by$ exists and satisfies
\begin{eqnarray*}\label{lse1}
d(1-4\sqrt{\delta(\log p)/n})\sigma^2 \leq \|\bX(\widehat{\bbeta}-\bbeta)\|_2^2 \leq d(1+4\sqrt{\delta(\log p)/n})\sigma^2.
\end{eqnarray*}
\end{lemma}


\begin{lemma}\label{lem2} 
Suppose that the $n \times p$ random matrix $\bX \sim N(\bzero, \bI_n \otimes \bSigma)$, where the eigenvalues of $\bSigma$ are bounded within the interval $[1/L, L]$ for some constant $L\geq 1$. For a positive constant $\xi$ and a set $S \subset \{1,\ldots,p\}$, the cone invertibility factor and the sign-restricted cone invertibility factor are defined as 
\begin{align*}
  F_2^* (\xi,S) & = \inf \left\{s^{1/2}\|\bX^\top\bX \bu/n\|_{\infty}/ \|\bu\|_{2}: \bu \in \mathcal{G}(\xi,S) \right\} \ \text{and}\\
  F_2 (\xi,S) & = \inf \left\{s^{1/2}\|\bX^\top\bX \bu/n\|_{\infty}/ \|\bu\|_{2}: \bu \in \mathcal{G}_{-}(\xi,S) \right\}, 
\end{align*}
respectively, where $s=|S|$, the cone $\mathcal{G}(\xi,S)=\{\bu\in \mathbb{R}^{p-1}:\|\bu_{S^c}\|_1 \leq \xi \|\bu_{S}\|_1 \neq 0\}$, and the sign-restricted cone $\mathcal{G}_{-}(\xi,S)=\{\bu\in \mathbb{R}^{p-1}:\|\bu_{S^c}\|_1 \leq \xi \|\bu_{S}\|_1 \neq 0,u_j\bx_j^\top\bX\bu/n\leq 0, \forall j\not\in S\}$.
Then for $s=o(n/\log (p))$, there exist positive constants $c_1$, $c_2$, $c_3$ such that with probability at least $1- c_2\exp(-c_3n)$,
\begin{align*}
F_2(\xi,S)\geq F_2^*  (\xi,S)  \geq c_1.
\end{align*}
\end{lemma}
\subsection{Proof of Proposition ~\ref{pro1}}
We first prove the equivalence of the optimization problems. For any $j$, $1 \leq j \leq p$, by definition we have
\begin{align}\label{pingguo}
&\widehat{\bpi}^{(j)}_{-j} = \mathop{\arg\min}_{\bb_{-j}} \Big\{(2n)^{-1}\big\| \bx_j - \bX_{-j} \bb_{-j}\big\|_2^2 +  \lambda_j \sum_{k \in \widehat{S}^c\backslash  \{j\} } v_{jk} \cdot |b_k| \Big\}.
\end{align}
Since the coefficients in the set $\widehat{S}\backslash  \{j\}$ are not penalized, the Karush-Kuhn-Tucker (KKT) condition for $\widehat{\bpi}^{(j)}_{\widehat{S}\backslash  \{j\}}$ gives
\begin{align*}
-\bX_{\widehat{S}\backslash  \{j\}}^\top(\bx_j-\bX_{\widehat{S}\backslash  \{j\}}\widehat{\bpi}^{(j)}_{\widehat{S}\backslash  \{j\}}-\bX_{\widehat{S}^c\backslash  \{j\}}\widehat{\bpi}^{(j)}_{\widehat{S}^c\backslash  \{j\}})=\bzero,
\end{align*}
which entails
\begin{align}\label{chun1}
\widehat{\bpi}^{(j)}_{\widehat{S}\backslash  \{j\}}=(\bX_{\widehat{S}\backslash  \{j\}}^\top\bX_{\widehat{S}\backslash  \{j\}})^{-1}\bX_{\widehat{S}\backslash  \{j\}}^\top(\bx_j-\bX_{\widehat{S}^c\backslash  \{j\}}\widehat{\bpi}^{(j)}_{\widehat{S}^c\backslash  \{j\}}).
\end{align}

When $\bb^{(j)}_{\widehat{S}\backslash  \{j\}}=(\bX_{\widehat{S}\backslash  \{j\}}^\top\bX_{\widehat{S}\backslash  \{j\}})^{-1}\bX_{\widehat{S}\backslash  \{j\}}^\top(\bx_j-\bX_{\widehat{S}^c\backslash  \{j\}}\bb^{(j)}_{\widehat{S}^c\backslash  \{j\}})$, we have
\begin{align}\label{qixi2}
\bx_j - \bX_{-j} \bb^{(j)}_{-j}=(\bI_n-\mathbb{\bP}_{\widehat{S}\backslash \{j\}})(\bx_j-\bX_{\widehat{S}^c\backslash  \{j\}}\bb^{(j)}_{\widehat{S}^c\backslash  \{j\}})=\bpsi_j^{(j)} - \bpsi_{\widehat{S}^c\backslash  \{j\}}^{(j)} \bb_{\widehat{S}^c\backslash  \{j\}}^{(j)}.
\end{align}
Thus, it follows from (\ref{chun1}) and (\ref{qixi2}) that the optimization problem  (\ref{pingguo}) is  equivalent to
\begin{align*}
\begin{cases}
\widehat{\bpi}^{(j)}_{\widehat{S}^c\backslash  \{j\}} = \mathop{\arg\min}_{\bb_{\widehat{S}^c\backslash  \{j\}}} \Big\{(2n)^{-1}\big\| \bpsi_j^{(j)} - \bpsi_{\widehat{S}^c\backslash  \{j\}}^{(j)} \bb_{\widehat{S}^c\backslash  \{j\}}\big\|_2^2 + \lambda_j \sum_{k \in \widehat{S}^c\backslash  \{j\} } v_{jk} \cdot |b_k| \Big\}, \\
\widehat{\bpi}^{(j)}_{\widehat{S}\backslash  \{j\}}=(\bX_{\widehat{S}\backslash  \{j\}}^\top\bX_{\widehat{S}\backslash  \{j\}})^{-1}\bX_{\widehat{S}\backslash  \{j\}}^\top(\bx_j-\bX_{\widehat{S}^c\backslash  \{j\}}\widehat{\bpi}^{(j)}_{\widehat{S}^c\backslash  \{j\}}).
\end{cases}
\end{align*}
Therefore, $\widehat{\bpi}^{(j)}_{\widehat{S}^c\backslash  \{j\}}$ and $\widehat{\bomega}_j$ in (\ref{penl}) are the optimal solutions of the same optimization problem. It yields that $\widehat{\bpi}^{(j)}_{\widehat{S}^c\backslash  \{j\}}$ and $\widehat{\bomega}_j$ are equivalent. Moreover, in view of (\ref{chun1}), we also have
\begin{align*}
\bx_j-\bX_{-j}\widehat{\bpi}^{(j)}_{-j}=(\bI_n-\mathbb{\bP}_{\widehat{S}\backslash \{j\}})(\bx_j-\bX_{\widehat{S}^c\backslash  \{j\}}\widehat{\bpi}^{(j)}_{\widehat{S}^c\backslash  \{j\}})=\bpsi_j^{(j)} - \bpsi_{\widehat{S}^c\backslash  \{j\}}^{(j)} \widehat{\bpi}^{(j)}_{\widehat{S}^c\backslash  \{j\}},
\end{align*}
which entails that $\bx_j-\bX_{-j}\widehat{\bpi}^{(j)}_{-j}$ is equivalent to the hybrid orthogonalization vector $\bz_j = \bpsi_j^{(j)} - \bpsi_{\widehat{S}^c\backslash  \{j\}}^{(j)} \widehat{\bomega}_j$ defined in (\ref{penl}) due to the equivalence of $\widehat{\bpi}^{(j)}_{\widehat{S}^c\backslash  \{j\}}$ and $\widehat{\bomega}_j$.

Second, although the penalized regression problem of $\widehat{\bpi}^{(j)}_{\widehat{S}^c\backslash  \{j\}}$ is not strictly convex and may not have a unique minimizer, we continue to show that different minimizers would give the same value of $\bpsi_{\widehat{S}^c\backslash  \{j\}}^{(j)} \widehat{\bpi}^{(j)}_{\widehat{S}^c\backslash  \{j\}}$. It implies that the value of the hybrid orthogonalization vector $\bz_j$ is unique. In what follows, we will prove this by contradiction.

Suppose that there are two solutions $\widehat{\bpi}^{(1j)}_{\widehat{S}^c\backslash  \{j\}}$ and $\widehat{\bpi}^{(2j)}_{\widehat{S}^c\backslash  \{j\}}$ such that $\bpsi_{\widehat{S}^c\backslash  \{j\}}^{(j)} \widehat{\bpi}^{(1j)}_{\widehat{S}^c\backslash  \{j\}}\neq\bpsi_{\widehat{S}^c\backslash  \{j\}}^{(j)}\widehat{\bpi}^{(2j)}_{\widehat{S}^c\backslash  \{j\}}$.
Denote by $t^*$ the minimum value of the optimization problem obtained by either $\widehat{\bpi}^{(1j)}_{\widehat{S}^c\backslash  \{j\}}$ or $\widehat{\bpi}^{(2j)}_{\widehat{S}^c\backslash  \{j\}}$. Then for any $0 < \alpha < 1$, the optimization problem attains a value of $t$ at $\alpha\widehat{\bpi}^{(1j)}_{\widehat{S}^c\backslash  \{j\}}+(1-\alpha)\widehat{\bpi}^{(2j)}_{\widehat{S}^c\backslash  \{j\}}$ as follows,
\begin{align*}
t= \frac{\big\|\bpsi_j^{(j)} - \bpsi_{\widehat{S}^c\backslash  \{j\}}^{(j)} (\alpha\widehat{\bpi}^{(1j)}_{\widehat{S}^c\backslash  \{j\}}+(1-\alpha)\widehat{\bpi}^{(2j)}_{\widehat{S}^c\backslash  \{j\}})\big\|_2^2}{2n} + \lambda_j \sum_{k \in \widehat{S}^c\backslash  \{j\} } v_{jk}  |\alpha\widehat{\pi}^{(1j)}_k+(1-\alpha)\widehat{\pi}^{(2j)}_k|,
\end{align*}
where $\widehat{\pi}_k^{(1j)}$ and $\widehat{\pi}_k^{(2j)}$  denote the $k$th elements of $\widehat{\bpi}^{(1j)}$ and $\widehat{\bpi}^{(2j)}$, respectively.

On one hand, due to the strict convexity of the function $f(\bx)=\|\by-\bx\|_2^2$, we get
\begin{align*}
& \big\|\bpsi_j^{(j)} - \bpsi_{\widehat{S}^c\backslash  \{j\}}^{(j)} (\alpha\widehat{\bpi}^{(1j)}_{\widehat{S}^c\backslash  \{j\}}+(1-\alpha)\widehat{\bpi}^{(2j)}_{\widehat{S}^c\backslash  \{j\}})\big\|_2^2\\
& < \alpha\big\|\bpsi_j^{(j)} - \bpsi_{\widehat{S}^c\backslash  \{j\}}^{(j)} \widehat{\bpi}^{(1j)}_{\widehat{S}^c\backslash  \{j\}}\big\|_2^2+(1-\alpha)\big\|\bpsi_j^{(j)} - \bpsi_{\widehat{S}^c\backslash  \{j\}}^{(j)}\widehat{\bpi}^{(2j)}_{\widehat{S}^c\backslash  \{j\}}\big\|_2^2.
\end{align*}
On the other hand, it follows from the convexity of the function $f(x)=|x|$ that
\begin{align*}
|\alpha\widehat{\pi}^{(1j)}_k+(1-\alpha)\widehat{\pi}^{(2j)}_k|\leq \alpha |\widehat{\pi}^{(1j)}_k|+(1-\alpha)|\widehat{\pi}^{(2j)}_k|.
\end{align*}
Thus, combining these two inequalities gives
$$t <\alpha t^* + (1-\alpha)t^*=t^*,$$
which means that $\alpha\widehat{\bpi}^{(1j)}_{\widehat{S}^c\backslash  \{j\}}+(1-\alpha)\widehat{\bpi}^{(2j)}_{\widehat{S}^c\backslash  \{j\}}$ attains a smaller value than $t^*$. This contradicts the definition of $t^*$, which completes the proof.

\subsection{Proof of Theorem \ref{haolei}}
We first prove the first part of this theorem. With some simple algebra, we have
\begin{eqnarray*}
|\Delta_j|=\Big|\sum_{k \in \widehat{S}^c\backslash  \{j\}}\frac{\bz_j^\top\bx_k\beta_k}{\bz_j^\top\bx_j}\Big|\leq \left(\max_{k \in \widehat{S}^c\backslash  \{j\}}\Big|\frac{\bz_j^\top\bx_k}{\bz_j^\top\bx_j}\Big|\right)\sum_{k \in {\widehat{S}^c}}|\beta_k|.
\end{eqnarray*}
In addition, by the KKT condition and the equality \eqref{zxb}, for any $k\in \widehat{S}^c\backslash  \{j\}$, we can get
\[|\bz_j^\top\bx_k|=|\bz_j^\top\bpsi_k^{(j)}| = | (\bpsi_j^{(j)} - \bpsi_{\widehat{S}^c\backslash  \{j\}}\widehat{\bomega}_j)^\top\bpsi_k^{(j)}| \leq \sqrt{n}\lambda_j \|\bpsi_k^{(j)}\|_2,\]
which along with  $\widehat{S}^c\subset {S}_1^c$ yields
\begin{align*}
|\Delta_j|\leq \left(\max_{k \in \widehat{S}^c\backslash  \{j\}}\Big|\frac{\sqrt{n}\lambda_j \|\bpsi_k^{(j)}\|_2}{\bz_j^\top\bx_j}\Big|\right)\sum_{k \in {{S}_1^c}}|\beta_k|.
\end{align*}

We proceed to prove the second part of this theorem.
With the aid of the  Lyapunov's condition \eqref{lya},  applying the Lyapunov's Central Limit Theorem gives
  \begin{align*}
\lim_{n\rightarrow\infty}\mathbb{P}(\frac{\bz_j^\top\bveps}{\sigma \|\bz_j\|_2}\leq  x)=\Phi(x),
\end{align*}
 where $x$ is any real number.
In view of $W_j=\frac{\bz_j^\top\bveps}{\bz_j^\top\bx_j}$, the above equality can be rewritten as
\begin{align*}
\lim_{n\rightarrow\infty}\mathbb{P}(W_j\leq \tau_j\sigma x)=\Phi(x),
\end{align*}
which completes the proof of this theorem.


\subsection{Proof of Corollary \ref{sanxia}}

The proof of this corollary is the same as that of Theorem \ref{haolei}, and, therefore, has been omitted.

\subsection{Proof of Proposition ~\ref{pro2}}

Since the design matrix $\bX$ here is random instead of a fixed one as that in \cite{Zhang2014}, the statistics related to $\bX$ are also random variables. We will first analyze the properties of some key statistics before proving Proposition ~\ref{pro2}.

\medskip
\noindent \emph{Part 1: Deviation bounds of $\|\bx_k\|_2$}. Under Condition \ref{fenbuz}, since $\|\bx_k\|_2^2/\sigma_{kk} \sim \chi_{(n)}^2$ for any $1 \leq k \leq p$ with $\sigma_{kk}$ denoting the $(k,k)$th entry of $\bSigma$, applying the following tail probability bound with $t = 4\sqrt{\delta \log (p)/n}$ for the chi-squared distribution with $n$ degrees of freedom \cite[inequality (93)]{RenZhou15}:
\begin{equation}\label{chi}
\mathbb{P}\big\{|\frac{\chi_{(n)}^2}{n} - 1| \geq t\big\} \leq 2 \exp\big(-nt(t \wedge 1)/8\big)
\end{equation}
gives that
\[[1 - 4\sqrt{\delta \log (p)/n}]\sigma_{kk}\leq \|\bx_k\|_2^2/n \leq [1 + 4\sqrt{\delta \log (p)/n}]\sigma_{kk}\]
holds with probability at least $1 - 2p^{-2\delta}$. By Condition \ref{fenbuz} that the eigenvalues of $\bPhi$ are within the interval $[1/L, L]$, with the aid of Lemma \ref{sigma}, we have $1/L\leq \sigma_{kk} \leq L$ for any $1 \leq k \leq p$. In view of $\log p=o(n)$ by Condition \ref{weidu}, it follows that for sufficiently large $n$, with probability at least $1 - 2p^{-2\delta}$,
\begin{equation}\label{changdu}
\widetilde{M}^{\ast}_1\leq\sqrt{[1 - 4\sqrt{\delta \log (p)/n}]/L}\leq\|\bx_k\|_2/\sqrt{n} \leq \sqrt{[1 + 4\sqrt{\delta \log (p)/n}]L} \leq \widetilde{M}_1,
\end{equation}
where $\widetilde{M}_1\geq \widetilde{M}^{\ast}_1$  are two positive constants.

Therefore, we have
\begin{equation*}
\begin{aligned}
&\mathbb{P}(\max_{k \neq j} \|\bx_k\|_2/\sqrt{n} > \widetilde{M}_1)\leq \sum_{k \neq j} \mathbb{P}(\|\bx_k\|_2/\sqrt{n} > \widetilde{M}_1)\leq p \cdot 2p^{-2\delta} = 2p^{1 - 2\delta},
\end{aligned}
\end{equation*}
as well as
\begin{equation*}
\begin{aligned}
&\mathbb{P}(\min_{k \neq j} \|\bx_k\|_2/\sqrt{n} < \widetilde{M}^{\ast}_1)\leq \sum_{k \neq j} \mathbb{P}(\|\bx_k\|_2/\sqrt{n} < \widetilde{M}^{\ast}_1)\leq p \cdot 2p^{-2\delta} = 2p^{1 - 2\delta}.
\end{aligned}
\end{equation*}
Combining the above two inequalities further yields
\begin{equation}\label{ycj}
\begin{aligned}
&\mathbb{P}\left(\bigcap_{k \neq j}\{\widetilde{M}^{\ast}_1\leq \|\bx_k\|_2/\sqrt{n} \leq \widetilde{M}_1\}\right)\geq 1-2p^{1 - 2\delta}-2p^{1 - 2\delta} =1- o(p^{1 - \delta}).
\end{aligned}
\end{equation}

\medskip

\noindent \emph{Part 2: Deviation bounds of $\|\bpsi_k^{(j)}\|_2$ for $k\in \widehat{S}^c$}. Assume that the index $k$ is the $k'$th index in the set $\widehat{S}^{(k)}=\widehat{S}\cup \{k\}$.
On one hand, in view of $\mathbb{\bP}_{\widehat{S}}=\bX_{\widehat{S}}(\bX_{\widehat{S}}^\top \bX_{\widehat{S}})^{-1}\bX_{\widehat{S}}^\top$, we
can rewrite $\bpsi_k^{(j)}$ as 
\begin{eqnarray}\label{guji1}
\bpsi_k^{(j)}=(\bI_n-\mathbb{\bP}_{\widehat{S}})\bx_k=\bx_k-\bX_{\widehat{S}}(\bX_{\widehat{S}}^\top \bX_{\widehat{S}})^{-1}\bX_{\widehat{S}}^\top\bx_k\triangleq\bx_k-\bX_{\widehat{S}} \widehat{\bgamma}_k^{(j)},
\end{eqnarray}
where $\widehat{\bgamma}_k^{(j)}=(\bX_{\widehat{S}}^\top \bX_{\widehat{S}})^{-1}\bX_{\widehat{S}}^\top\bx_k$ is the least squares estimate of the coefficient vector of regressing $\bx_k$ on $\bX_{\widehat{S}}$.

On the other hand, since the rows of $\bX_{\widehat{S}^{(k)}}$ are independent and normally distributed with $N(0, (\bTheta^{(k)})^{-1})$ for  $\bTheta^{(k)}=(\bSigma_{\widehat{S}^{(k)},\widehat{S}^{(k)}})^{-1}=(\theta^{(k)}_{ij})_{(d+1)\times (d+1)}$, the conditional distribution of $\bx_{k}$ given $\bX_{\widehat{S}}$ follows the Gaussian distribution
$$\bx_{k}|\bX_{\widehat{S}}\sim N(\bX_{\widehat{S}}\bgamma_k^{(j)}, (1/\theta^{(k)}_{k'k'})\bI_n)$$
with $\bgamma_k^{(j)}=-(\bTheta^{(k)}_{k',-k'})^\top/\theta^{(k)}_{k'k'}$. Here the $d$-dimensional coefficient vector $\bTheta^{(k)}_{k',-k'}$ is the $k'$th row of $\bTheta^{(k)}$ with the $k'$th component removed.
Based on this conditional distribution, we know that the random error vector
\begin{eqnarray}\label{fenbu}
\brho_{k}^{(j)} = \bx_{k}-\bX_{\widehat{S}}\bgamma_k^{(j)}
\end{eqnarray}
is independent of $\bX_{\widehat{S}}$ and  follows the distribution $N(0,(1/\theta^{(k)}_{k'k'})\bI_n)$.
It is worth mentioning that $\theta^{(k)}_{k'k'}\in [1/L, L]$ can be easily deduced from Condition \ref{fenbuz} and Lemma \ref{sigma}.

Now we provide the deviation of $\|\bpsi_k^{(j)}\|_2$ from its population counterpart $\|\brho_{k}^{(j)}\|_2$.
In view of  Lemma \ref{lse-lem},  combining inequalities (\ref{guji1}) and (\ref{fenbu}) above yields
\begin{eqnarray*}
&\mathbb{P}&\left\{\left(\frac{d(1-4\sqrt{\delta(\log p)/n})}{\theta^{(k)}_{k'k'}}\right)^{1/2} \leq \|\bX_{\widehat{S}}(\widehat{\bgamma}_k^{(j)}-\bgamma_k^{(j)})\|_2
\leq \left(\frac{d(1+4\sqrt{\delta(\log p)/n})}{\theta^{(k)}_{k'k'}}\right)^{1/2}\right\}\\
&\ge& 1- 2p^{-2\delta}.
\end{eqnarray*}
Moreover, in view of (\ref{guji1}) and (\ref{fenbu}), we also have
$$\brho_k^{(j)}-\bpsi_k^{(j)}=\bX_{\widehat{S}}(\widehat{\bgamma}_k^{(j)}-\bgamma_k^{(j)}),$$
which together with the above inequality gives
\begin{eqnarray}\nonumber
&&\mathbb{P}\bigg\{\bigg(\frac{d(1-4\sqrt{\delta(\log p)/n})}{\theta^{(k)}_{k'k'}}\bigg)^{1/2} \leq \|\brho_k^{(j)}-\bpsi_k^{(j)}\|_2\leq \bigg(\frac{d(1+4\sqrt{\delta(\log p)/n})}{\theta^{(k)}_{k'k'}}\bigg)^{1/2}\bigg\}\\ \label{zhaohui}
&\ge& 1- 2p^{-2\delta}.
\end{eqnarray}
Since $\big| \|\bpsi_k^{(j)}\|_2 - \|\brho_k^{(j)}\|_2\big| \leq \|\brho_k^{(j)}-\bpsi_k^{(j)}\|_2$, the above inequality entails
\begin{eqnarray}\label{gujifen}
\mathbb{P} \left\{\big|\|\bpsi_k^{(j)}\|_2 - \|\brho_k^{(j)}\|_2\big|  \ge \left(\frac{d(1+4\sqrt{\delta(\log p)/n})}{\theta^{(k)}_{k'k'}}\right)^{1/2}\right\}\le 2p^{-2\delta}.
\end{eqnarray}

In order to proceed, we will then construct the deviation bounds of $\|\brho_k^{(j)}\|_2$. Note that $\|\brho_k^{(j)}\|_2^2/(1/\theta^{(k)}_{k'k'}) \sim \chi_{(n)}^2$. By applying (\ref{chi}) with $t = 4\sqrt{\delta \log (p)/n}$ for the chi-squared distribution with $n$ degrees of freedom, we  have
\begin{equation*}
\left(\frac{n(1 - 4\sqrt{\delta\log (p)/n})}{\theta^{(k)}_{k'k'}}\right)^{1/2} \leq \|\brho_k^{(j)}\|_2 \leq \left(\frac{n(1 + 4\sqrt{\delta\log (p)/n})}{\theta^{(k)}_{k'k'}}\right)^{1/2}
\end{equation*}
holds with probability at least $1 - 2p^{-2\delta}$.  This inequality together with (\ref{gujifen}) entails that with probability at least $1 - 4p^{-2\delta}$,
\begin{eqnarray*}
&&\frac{(n(1 - 4\sqrt{\delta\log (p)/n}))^{1/2}-(d(1+4\sqrt{\delta(\log p)/n}))^{1/2}}{(\theta^{(k)}_{k'k'})^{1/2}}\\
&\leq &\|\bpsi_k^{(j)}\|_2 \leq \frac{(n^{1/2}+d^{1/2})(1 + 4\sqrt{\delta\log (p)/n})^{1/2}}{(\theta^{(k)}_{k'k'})^{1/2}}.
\end{eqnarray*}

Since $\log (p)=o(n)$ and $d=o(n/\log (p))$, which entails $d=o(n)$,
it follows that for sufficiently large $n$, with probability at least $1 - 4p^{-2\delta}$,
\begin{equation}\label{changdu1}
\widetilde{M}^{\ast}_2\leq\|\bpsi_k^{(j)}\|_2/\sqrt{n} \leq \widetilde{M}_2,
\end{equation}
where $\widetilde{M}_2\geq \widetilde{M}^{\ast}_2$ are two positive constants.
Thus, we have
\begin{equation*}
\begin{aligned}
&\mathbb{P}(\max_{k\in\widehat{S}^c\backslash  \{j\}} \|\bpsi_k^{(j)}\|_2/\sqrt{n} > \widetilde{M}_2)\leq \sum_{k \neq j} \mathbb{P}(\|\bpsi_k^{(j)}\|_2/\sqrt{n} > \widetilde{M}_2)\leq p \cdot 4p^{-2\delta} = 4p^{1 - 2\delta}
\end{aligned}
\end{equation*}
and
\begin{equation*}
\begin{aligned}
&\mathbb{P}(\min_{k\in\widehat{S}^c\backslash  \{j\}} \|\bpsi_k^{(j)}\|_2/\sqrt{n} < \widetilde{M}^{\ast}_2)\leq \sum_{k \neq j} \mathbb{P}(\|\bpsi_k^{(j)}\|_2/\sqrt{n} < \widetilde{M}^{\ast}_2)\leq p \cdot 4p^{-2\delta} = 4p^{1 - 2\delta},
\end{aligned}
\end{equation*}
respectively. Combining the above two inequalities further yields
\begin{equation}\label{ycj1}
\begin{aligned}
&\mathbb{P}\left(\bigcap_{k\in\widehat{S}^c\backslash  \{j\}}\{\widetilde{M}^{\ast}_2\leq \|\bpsi_k^{(j)}\|_2/\sqrt{n} \leq \widetilde{M}_2\}\right)\geq 1-4p^{1 - 2\delta}-4p^{1 - 2\delta} =1- o(p^{1 - \delta}).
\end{aligned}
\end{equation}

\medskip
We proceed to  prove Proposition ~\ref{pro2}.
Since the following argument applies to any $j$, $1 \leq j \leq p$, we drop the superscript $(j)$ in $F_2^{(j)} (\xi,S_{j*})$ and $\mathcal{G}_{-}^{(j)}(\xi,S_{j*})$ for notational simplicity. Recall that the sign-restricted cone invertibility factor here is defined as
\begin{align*}
  F_2 (\xi,S_{j*})  = \inf \left\{|S_{j*}|^{1/2}\|(\bpsi_{\widehat{S}^c\backslash  \{j\}}^{(j)})^\top\bpsi_{\widehat{S}^c\backslash  \{j\}}^{(j)} \bu/n\|_{\infty}/ \|\bu\|_{2}: \bu \in \mathcal{G}_{-}(\xi,S_{j*}) \right\} 
\end{align*}
with the sign-restricted cone 
$$\mathcal{G}_{-}(\xi,S_{j*})=\{\bu: \|\bu_{S_{j*}^c}\|_1 \leq \xi \|\bu_{S_{j*}}\|_1, \ u_k(\bpsi_k^{(j)})^\top\bpsi_{\widehat{S}^c\backslash  \{j\}}^{(j)} \bu/n\leq 0, \ \forall \ k \in S_{j*}^c\}.$$
In what follows, we mainly present the proof when $j \in \widehat{S}^c$ and the same argument applies to $j \in \widehat{S}$.  

First of all, by the inequality \eqref{zhaohui}, we have
\begin{eqnarray*}
&&\mathbb{P}\left\{\left(\frac{d(1-4\sqrt{\delta(\log p)/n})}{\theta^{(k)}_{k'k'}}\right)^{1/2} \leq \|\brho_k^{(j)}-\bpsi_k^{(j)}\|_2\leq \left(\frac{d(1+4\sqrt{\delta(\log p)/n})}{\theta^{(k)}_{k'k'}}\right)^{1/2}\right\}\\ 
&\ge& 1- 2p^{-2\delta}.
\end{eqnarray*}
Since $\log (p)=o(n)$, it follows that for sufficiently large $n$, with probability at least $1 - 2p^{-2\delta}$,
\begin{equation*}
\|\brho_k^{(j)}-\bpsi_k^{(j)}\|_2/\sqrt{d} \leq C,
\end{equation*}
where $C$ is some positive constant.
Thus, we further have
\begin{equation}\label{sunl}
\begin{aligned}
&\mathbb{P}(\max_{k\in\widehat{S}^c\backslash  \{j\}} \|\brho_k^{(j)}-\bpsi_k^{(j)}\|_2/\sqrt{d} > C)\leq \sum_{k \neq j} \mathbb{P}(\|\brho_k^{(j)}-\bpsi_k^{(j)}\|_2/\sqrt{d} > C)\leq p \cdot 2p^{-2\delta} = 2p^{1 - 2\delta}.
\end{aligned}
\end{equation}

We then prove the conclusions of this proposition. For any $\bu\in \mathbb{R}^{p-d-1}$, we will first give the deviation of $\|\bpsi_k^{(j)}\bu\|_2$ from its population counterpart $\|\brho_{k}^{(j)}\bu\|_2$.
For any matrix $\bH=(\bh_1,\ldots,\bh_p)\in \mathbb{R}^{n\times p}$ and any vector $\bt=(t_1, \ldots, t_p)^\top\in \mathbb{R}^{p}$, it is clear that
$$\|\bH\bt\|_2=\|t_1\bh_1+\ldots+t_p\bh_p\|_2\leq |t_1|\|\bh_1\|_2+\ldots+|t_p|\|\bh_p\|_2\leq (\max_{1\leq i\leq p}\|\bh_i\|_2)\|\bt\|_1.$$
Based on this inequality, we have
$$\|\brho^{(j)}_{\widehat{S}^c\backslash  \{j\}}\bu\|_2-\|\bpsi_{\widehat{S}^c\backslash  \{j\}}^{(j)}\bu\|_2\leq\|(\brho^{(j)}_{\widehat{S}^c\backslash  \{j\}}-\bpsi_{\widehat{S}^c\backslash  \{j\}}^{(j)})\bu\|_2\leq \max_{k\in \widehat{S}^c\backslash  \{j\}}\|\brho_k^{(j)}-\bpsi_k^{(j)}\|_2\|\bu\|_1.$$
Therefore, for any  $\bu\in \mathcal{G}(\xi,S_{j*})$ with $\mathcal{G}(\xi,S_{j*})=\{\bu\in \mathbb{R}^{p-d-1}:\|\bu_{S_{j*}^c}\|_1 \leq \xi \|\bu_{S_{j*}}\|_1\}$,  in view of \eqref{sunl}, $|S_{j*}| \leq s_0$, and $\|\bu_{S_{j*}}\|_1\leq s_0^{1/2}\|\bu_{S_{j*}}\|_2$, with probability at least $1 - o(p^{1-\delta})$,
\begin{align}\nonumber
&\|\bpsi_{\widehat{S}^c\backslash  \{j\}}^{(j)}\bu\|_2\geq \|\brho^{(j)}_{\widehat{S}^c\backslash  \{j\}}\bu\|_2-C d^{1/2}\|\bu\|_1\\\label{tieg}
\geq \ & \|\brho^{(j)}_{\widehat{S}^c\backslash  \{j\}}\bu\|_2-C d^{1/2} (1+\xi) \|\bu_{S_{j*}}\|_1
\geq \|\brho^{(j)}_{\widehat{S}^c\backslash  \{j\}}\bu\|_2- C(1+\xi)(s_0d)^{1/2} \|\bu_{S_{j*}}\|_2.
\end{align}

Now we turn our attention to $\brho_{\widehat{S}^c\backslash  \{j\}}^{(j)}$. By Lemma \ref{tongfb}, we have
$$\brho^{(j)}_{\widehat{S}^c\backslash  \{j\}}\sim N(0, \bI_n \otimes [\bSigma_{\widehat{S}^c\backslash  \{j\}\widehat{S}^c\backslash  \{j\}}-\bSigma_{\widehat{S}^c\backslash  \{j\}\widehat{S}}\bSigma_{\widehat{S}\widehat{S}}^{-1}\bSigma_{\widehat{S}\widehat{S}^c\backslash  \{j\}}])$$
with $\bSigma_{\widehat{S}^c\backslash  \{j\}\widehat{S}^c\backslash  \{j\}}-\bSigma_{\widehat{S}^c\backslash  \{j\}\widehat{S}}\bSigma_{\widehat{S}\widehat{S}}^{-1}\bSigma_{\widehat{S}\widehat{S}^c\backslash  \{j\}}=\bPhi_{\widehat{S}^c\backslash  \{j\}\widehat{S}^c\backslash  \{j\}}^{-1}$.
Since the eigenvalues of $\bPhi$ are bounded within the interval $[1/L, L]$ for some constant $L\geq 1$, we can deduce that the eigenvalues of $\bPhi_{\widehat{S}^c\backslash  \{j\}\widehat{S}^c\backslash  \{j\}}^{-1}$ are also bounded within the interval $[1/L, L]$ by Lemma \ref{sigma}.
Thus, for any $\bu\in \mathcal{G}(\xi,S_{j*})$, applying the same argument as that of \eqref{wuyule} in Lemma \ref{lem2} yields that
\begin{align*}
\frac{\|\brho^{(j)}_{\widehat{S}^c\backslash  \{j\}}\bu\|_2}{n^{1/2}}\geq \left\{\frac{1}{4L^{1/2}}-9(1+\xi)L^{1/2}\sqrt{\frac{s_0 \log (p)}{n}}\right\}\|\bu\|_2
\end{align*}
holds with probability at least $1- c_2^*\exp(-c_3^*n)$, where $c_2^*$ and $c_3^*$ are some positive constants.

Together with (\ref{tieg}), it follows that for any $\bu\in \mathcal{G}(\xi,S_{j*})$,
\begin{align}\nonumber
\frac{\|\bpsi_{\widehat{S}^c\backslash  \{j\}}^{(j)}\bu\|_2}{n^{1/2}} &\geq \frac{\|\brho^{(j)}_{\widehat{S}^c\backslash  \{j\}}\bu\|_2}{n^{1/2}}- C(1+\xi)(\frac{s_0d}{n})^{1/2} \|\bu_{S_{j*}}\|_2\\ \label{hhgai}
&\geq \left\{\frac{1}{4L^{1/2}}-9(1+\xi)L^{1/2}\sqrt{\frac{s_0 \log (p)}{n}}-C(1+\xi)(\frac{s_0d}{n})^{1/2}\right\}\|\bu\|_2
\end{align}
holds with probability at least $1- c_2^*\exp(-c_3^*n) - o(p^{1-\delta})$.
Moreover, since $\log(p)/n=o(1)$, it entails that
$$\frac{\exp(-c_3^*n)}{p^{1-\delta}}=\frac{p^{\delta-1}}{\exp(c_3^*n)}=\exp\{(\delta-1)\log p -c_3^*n\}=o(1),$$
which implies that $c_2^*\exp(-c_3^*n)=o(p^{1-\delta})$.

Therefore, based on this fact, together with the assumptions of $s_0 (\log p)/n = o(1)$ and $d = o(n/s_0)$,  it follows from \eqref{hhgai} that
\begin{align}\label{hhgb}
RE_2 (\xi,S_{j*})  \geq (c^{*})^{1/2}
\end{align}
holds with probability at least $1- o(p^{1-\delta})$, where $c^{*}$ is some constant and $RE_2 (\xi,S_{j*})$ is defined as
$$RE_2 (\xi,S_{j*})= \inf \left\{\frac{\|\bpsi_{\widehat{S}^c\backslash  \{j\}}^{(j)}\bu\|_{2}}{n^{1/2} \|\bu\|_{2}}: \bu \in \mathcal{G}(\xi,S_{j*}) \right\}.$$
Note that the cone $\mathcal{G}(\xi,S_{j*})=\{\bu\in \mathbb{R}^{p-d-1}:\|\bu_{S_{j*}^c}\|_1 \leq \xi \|\bu_{S_{j*}}\|_1\}$.

Finally, similar to the proof of Lemma \ref{lem2}, by \cite[Proposition 5]{ye2010}, we can get
\begin{align*}
F_2 (\xi,S_{j*})\geq F_2^* (\xi,S_{j*}) \geq RE_2^2 (\xi,S_{j*}) \geq c^{*}
\end{align*}
holds with probability at least $1- o(p^{1-\delta})$,
where the cone invertibility factor $F_2^* (\xi,S_{j*})$ is defined as 
\begin{align*}
  F_2^* (\xi,S_{j*})  = \inf \left\{|S_{j*}|^{1/2}\|(\bpsi_{\widehat{S}^c\backslash  \{j\}}^{(j)})^\top\bpsi_{\widehat{S}^c\backslash  \{j\}}^{(j)} \bu/n\|_{\infty}/ \|\bu\|_{2}: \bu \in \mathcal{G}(\xi,S_{j*}) \right\}. 
\end{align*}
It completes the proof of Proposition ~\ref{pro2}.

\subsection{Proof of Theorem \ref{theob}}

Under the conditions of Theorem  \ref{theob}, it is clear that the results in \emph{Part 1} and \emph{Part 2} of the proof of Proposition \ref{pro2} are still valid.
In what follows, we mainly present the proof when $j \in \widehat{S}^c$ and the same argument applies to $j \in \widehat{S}$ for a given set $\widehat{S}$. Throughout the proof,  $\lambda_j=(1+\epsilon)\sqrt{(2\delta/(n\phi_{jj})) \log(p)}$,  where the constant $\epsilon=\widetilde{M}(\xi'+1)/\{\widetilde{M}^{\ast}(\xi'-1)\} - 1>0$ with two positive constants $\widetilde{M}\geq \widetilde{M}^{\ast}$, and $\phi_{jj}$ is defined in Condition \ref{fenbuz}. In addition, $C$ is some positive constant which can vary from expression to expression.

\medskip
\noindent\emph{Part 1: Deviation bounds of $\|\widehat{\bomega}_j-{\bomega}_j\|_2$.}
Let $\widetilde{M}=\max(\widetilde{M}_1,\widetilde{M}_2)$ and $\widetilde{M}^{\ast}=\min(\widetilde{M}^{\ast}_1,\widetilde{M}^{\ast}_2)$, respectively. Recall that we have the conditional distribution
$$\bx_j=\bX_{\widehat{S}}{\bpi}^{(j)}_{\widehat{S}}+\bX_{\widehat{S}^c\backslash  \{j\}}{\bpi}^{(j)}_{\widehat{S}^c\backslash  \{j\}}+\textbf{e}_j,$$
where $\textbf{e}_j=\bx_j- \mathbb{E}(\bx_j|\bX_{-j})$  is independent of $\bX_{-j}$ and follows the distribution $N(\bzero,(1/\phi_{jj})\bI_n)$.
Thus, it follows that
$$(\bI_n-\mathbb{\bP}_{\widehat{S}})\bx_j=(\bI_n-\mathbb{\bP}_{\widehat{S}})(\bX_{\widehat{S}}{\bpi}^{(j)}_{\widehat{S}}+\bX_{\widehat{S}^c\backslash  \{j\}}{\bpi}^{(j)}_{\widehat{S}^c\backslash  \{j\}}+\textbf{e}_j),$$
which can be simplified as
\begin{align}\label{lud}
\bpsi_j^{(j)}=\bpsi_{\widehat{S}^c\backslash  \{j\}}^{(j)} \bomega_j+\textbf{e}_j^*
\end{align}
according to the definition of $\bpsi_k^{(j)}$. Here $\bomega_j={\bpi}^{(j)}_{\widehat{S}^c\backslash  \{j\}}=\bPhi_{j,{\widehat{S}^c\backslash  \{j\}}}^\top/\phi_{jj}$ and $\textbf{e}_j^*=(\bI_n-\mathbb{\bP}_{\widehat{S}})\textbf{e}_j$.

Denote by $S_{j*}=\{1 \leq k \leq p: k\in {\widehat{S}^c\backslash  \{j\}}, \phi_{jk}\neq 0\}$, $S_{j*}^c=\widehat{S}^c\backslash  (j\cup S_{j*})$, and $z_{\infty}^{j*}=\|(\bpsi_{\widehat{S}^c\backslash  \{j\}}^{(j)})^\top(\bpsi_j^{(j)}-\bpsi_{\widehat{S}^c\backslash  \{j\}}^{(j)}\bomega_j)/n\|_{\infty}$. In what follows, we  first use similar arguments as those in the proof of Theorem \ref{theoa} to construct bounds for $\|\bpsi_{\widehat{S}^c\backslash \{j\}}^{(j)}(\widehat{\bomega}_j-{\bomega}_j)\|_2$ and the main difficulty here is that the predictor matrix $\bpsi_{\widehat{S}^c\backslash \{j\}}^{(j)}$ is much more complex than a usual design matrix with i.i.d. rows.
First of all, by the same argument as that in (\ref{yuebing}) gives
\begin{align*}
\bv^\top(\bpsi_{\widehat{S}^c\backslash  \{j\}}^{(j)})^\top\bpsi_{\widehat{S}^c\backslash  \{j\}}^{(j)}(\widehat{\bomega}_j-{\bomega}_j)/n\leq (z_{\infty}^{j*}+ \widetilde{M}\lambda_j)\|\bv_{J_j}\|_1+(z_{\infty}^{j*}- \widetilde{M}^{\ast}\lambda_j)\|\bv_{J_j^c}\|_1
\end{align*}
for all vectors $\bv \in \mathbb{R}^{p-d-1}$ satisfying $\text{sign}(\bv_{S_{j*}^c})=\text{sign}((\widehat{\bomega}_j-{\bomega}_j)_{S_{j*}^c})=\text{sign}((\widehat{\bomega}_j)_{S_{j*}^c})$.
Based on the above inequality, if  $z_{\infty}^{j*}<\widetilde{M}^{\ast}\lambda_j$, it follows from $\bv=\widehat{\bomega}_j-{\bomega}_j$ that
$$(z_{\infty}^{j*}+ \widetilde{M}\lambda_j)\|(\widehat{\bomega}_j-{\bomega}_j)_{S_{j*}}\|_1+(z_{\infty}^{j*}- \widetilde{M}^{\ast}\lambda_j)\|(\widehat{\bomega}_j-{\bomega}_j)_{S_{j*}^c}\|_1\geq0,$$
which entails that
\begin{align}\label{zixidian1}
\|(\widehat{\bomega}_j-{\bomega}_j)_{S_{j*}^c}\|_1 \leq \xi \|(\widehat{\bomega}_j-{\bomega}_j)_{S_{j*}}\|_1
\end{align}
for  $\xi\geq (\widetilde{M}\lambda_{j}+z_{\infty}^{j*})/(\widetilde{M}^{\ast}\lambda_j-z_{\infty}^{j*})$.

In addition, if $z_{\infty}^{j*}\leq\widetilde{M}^{\ast}\lambda_j$, using the same argument as that in (\ref{henzixi}) gives
\begin{align}\label{henzixi1}
 (\widehat{\omega}_{j,k}-{\omega}_{j,k})(\bpsi_k^{(j)})^\top\bpsi_{\widehat{S}^c\backslash  \{j\}}^{(j)}(\widehat{\bomega}_j-{\bomega}_j)/n\leq 0
\end{align}
for any $ k\in J_j^c$, where  $\homega_{j,k}$ is the $k$th component of $\widehat{\bomega}_{j}$.
In view of \eqref{zixidian1} and \eqref{henzixi1}, if $z_{\infty}^{j*}<\widetilde{M}^{\ast}\lambda_j$, we can get
\begin{align}\label{chaozixi1}
\widehat{\bomega}_j-{\bomega}_j\in \mathcal{G}_{-}(\xi,S_{j*})
 \end{align}
 for $\xi\geq (\widetilde{M}\lambda_j+z_{\infty}^{j*})/(\widetilde{M}^{\ast}\lambda_j-z_{\infty}^{j*})$. Here the sign-restricted cone $\mathcal{G}_{-}(\xi,S_{j*})$ is defined as 
$$\mathcal{G}_{-}(\xi,S_{j*})=\{\bu\in \mathbb{R}^{p-d-1}:\|\bu_{S_{j*}^c}\|_1 \leq \xi \|\bu_{S_{j*}}\|_1, u_k(\bpsi_k^{(j)})^\top\bpsi_{\widehat{S}^c\backslash  \{j\}}^{(j)} \bu/n\leq 0, \forall k\in S_{j*}^c\}.$$

In order to proceed, our analysis will then be conditional on the event $\mathcal{E}^{'''*}= \{z_{\infty}^{j*}\leq \lambda_j\widetilde{M}^{\ast}(\xi'-1)/(\xi'+1)\}$ for some constant $\xi'>1$. Note that
$z_{\infty}^{j*}\leq \lambda_j\widetilde{M}^{\ast}(\xi'-1)/(\xi'+1)< \widetilde{M}^{\ast}\lambda_j$
holds on this event, which together with \eqref{chaozixi1} implies that
\begin{align*}
\widehat{\bomega}_j-{\bomega}_j\in \mathcal{G}_{-}(\xi,S_{j*})
 \end{align*}
for some constant $\xi\geq (\widetilde{M}\lambda_j+z_{\infty}^{j*})/(\widetilde{M}^{\ast}\lambda_j-z_{\infty}^{j*})$.
Similarly, we first clarify that the magnitude of $(\widetilde{M}\lambda_j+z_{\infty}^{j*})/(\widetilde{M}^{\ast}\lambda_j-z_{\infty}^{j*})$ is around a constant level. By the same argument
as that in \eqref{yanzi} gives
\begin{align*}
1\leq \frac{\widetilde{M}\lambda_j+z_{\infty}^{j*}}{\widetilde{M}^{\ast}\lambda_j-z_{\infty}^{j*}}\leq 1+\frac{\xi'-1}{2}+\frac{(\xi'+1)\widetilde{M}-2\widetilde{M}^{\ast}}{2\widetilde{M}^{\ast}}.
\end{align*}

Moreover, since $z_{\infty}^{j*}\leq \lambda_j\widetilde{M}^{\ast}(\xi'-1)/(\xi'+1)\leq \lambda_j\widetilde{M}(\xi'-1)/(\xi'+1)$, it follows that
\begin{align*}
\|(\bpsi_{\widehat{S}^c\backslash  \{j\}}^{(j)})^\top\bpsi_{\widehat{S}^c\backslash  \{j\}}^{(j)}(\widehat{\bomega}_j-{\bomega}_j)/n\|_{\infty}\leq \widetilde{M}\lambda_j+z_{\infty}^{j*}\leq 2\widetilde{M}\lambda_j\xi'/(\xi'+1).
\end{align*}
Thus, in view of $|S_{j*}|\leq s_0$, combining the above results gives
\begin{align}\label{xuezhuan1}
\|\widehat{\bomega}_j-{\bomega}_j\|_2\leq \frac{|S_{j*}|^{1/2}(\widetilde{M}\lambda_j+z_{\infty}^{j*})}{F_2(\xi,S_{j*})}\leq\frac{\{2\widetilde{M}\xi'/(\xi'+1)\}s_0^{1/2}\lambda_j}{F_2(\xi,S_{j*})},
\end{align}
where the sign-restricted cone invertibility factor $F_2(\xi,S_{j*})$ is defined as 
\begin{align*}
  F_2 (\xi,S_{j*})  = \inf \left\{|S_{j*}|^{1/2}\|(\bpsi_{\widehat{S}^c\backslash  \{j\}}^{(j)})^\top\bpsi_{\widehat{S}^c\backslash  \{j\}}^{(j)} \bu/n\|_{\infty}/ \|\bu\|_{2}: \bu \in \mathcal{G}_{-}(\xi,S_{j*}) \right\}. 
\end{align*}

We then derive the probability of the event $\mathcal{E}^{'''*}$.
Since $\textbf{e}_j$ is independent of $\bX_{\widehat{S}}$ and $\bX_{\widehat{S}^c\backslash  \{j\}}$, it entails that $\textbf{e}_j$ is independent of $\bpsi_{\widehat{S}^c\backslash  \{j\}}^{(j)}$.
In addition, it follows from the definition of $\textbf{e}_j^* =(\bI_n-\mathbb{\bP}_{\widehat{S}})\textbf{e}_j$ that
$$(\bpsi_{\widehat{S}^c\backslash  \{j\}}^{(j)})^\top\textbf{e}_j^*=(\bpsi_{\widehat{S}^c\backslash  \{j\}}^{(j)})^\top\textbf{e}_j.$$
Thus, together with $\|\bpsi_k^{(j)}\|_2/\sqrt{n} \leq \widetilde{M}$ for any $k \in \widehat{S}^c\backslash  \{j\}$ on the event $\mathcal{E}''$ and $\textbf{e}_j$ follows the distribution $N(0,(1/\phi_{jj})\bI_n)$,
when $\lambda_j=\widetilde{M}\sqrt{(2\delta/(n\phi_{jj})) \log(p)}(\xi'+1)/\{\widetilde{M}^{\ast}(\xi'-1)\}$, using the same argument as that in (\ref{luy1}) gives
\begin{align}\label{hhaha}
1-\mathbb{P}\left(\mathcal{E}^{'''*}\right)\leq o(p^{1-\delta}).
\end{align}

We now show the magnitude of $F_2(\xi,S_{j*})$.
By Proposition \ref{pro2}, under the assumption that $d=o(n/s_0)$, there is some positive constant $c^{*}$ such that
\begin{align*}
F_2 (\xi,S_{j*})  \geq c^{*}
\end{align*}
holds with probability at least $1- o(p^{1-\delta})$.
Moreover, with the aid of  (\ref{ycj1}) and \eqref{hhaha}, we get
\begin{align*}
&\mathbb{P}(\mathcal{E}''\cap\mathcal{E}^{'''*})=1-\mathbb{P}\left((\mathcal{E}'')^c\cup(\mathcal{E}^{'''*})^c\right)\geq  1-\mathbb{P}((\mathcal{E}'')^c)-\mathbb{P}((\mathcal{E}^{'''*})^c)=1-o(p^{1-\delta}).
\end{align*}
Thus, combining these results with \eqref{xuezhuan1} gives
\begin{eqnarray}\label{sbx}
\mathbb{P} \left\{\|\widehat{\bomega}_j-{\bomega}_j\|_2\leq Cs_0^{1/2}\lambda_j\right\}\geq 1-o(p^{1-\delta}).
\end{eqnarray}

\medskip

\noindent \emph{Part 2: Deviation bounds of $n^{-1/2}\|\bpsi_{\widehat{S}^c\backslash  \{j\}}^{(j)}(\widehat{\bomega}_j-{\bomega}_j)\|_2$ and $\|\widehat{\bomega}_j-{\bomega}_j\|_q$ for $q\in[1,2)$.}
We first construct bounds for $\|\bpsi_{\widehat{S}^c\backslash  \{j\}}^{(j)}(\widehat{\bomega}_j-{\bomega}_j)\|_2$. Some simple algebra shows that
\begin{align}\label{sdd}
\|\bpsi_{\widehat{S}^c\backslash  \{j\}}^{(j)}(\widehat{\bomega}_j-{\bomega}_j)\|_2^2\leq \|(\bpsi_{\widehat{S}^c\backslash  \{j\}}^{(j)})^\top\bpsi_{\widehat{S}^c\backslash  \{j\}}^{(j)}(\widehat{\bomega}_j-{\bomega}_j)\|_{\infty}\|(\widehat{\bomega}_j-{\bomega}_j)\|_1.
\end{align}
Conditional on the event $\mathcal{E}^{'''*}=\{z_{\infty}^{j*}\leq \lambda_j\widetilde{M}^{\ast}(\xi'-1)/(\xi'+1)\}$, it follows that
$$\|(\bpsi_{\widehat{S}^c\backslash  \{j\}}^{(j)})^\top\bpsi_{\widehat{S}^c\backslash  \{j\}}^{(j)}(\widehat{\bomega}_j-{\bomega}_j)/n\|_{\infty}\leq \widetilde{M}\lambda_j+z_{\infty}^{j*} \leq  \widetilde{M}\lambda_j+\lambda_j\widetilde{M}^{\ast}(\xi'-1)/(\xi'+1),$$
which further yields
\begin{align}\label{sddr}
\|(\bpsi_{\widehat{S}^c\backslash  \{j\}}^{(j)})^\top\bpsi_{\widehat{S}^c\backslash  \{j\}}^{(j)}(\widehat{\bomega}_j-{\bomega}_j)\|_{\infty}\leq n \{\widetilde{M}+\widetilde{M}^{\ast}(\xi'-1)/(\xi'+1)\}\lambda_j.
\end{align}

Moreover, since $\widehat{\bomega}_j-{\bomega}_j\in \mathcal{G}_{-}(\xi,S_{j*})$, we have
\begin{align*}
\|\widehat{\bomega}_j-{\bomega}_j\|_1&=\|(\widehat{\bomega}_j-{\bomega}_j)_{S_{j*}^c}\|_1 + \|(\widehat{\bomega}_j-{\bomega}_j)_{S_{j*}}\|_1\leq (1+\xi) \|(\widehat{\bomega}_j-{\bomega}_j)_{S_{j*}}\|_1\\
&\leq (1+\xi)s_0^{1/2} \|(\widehat{\bomega}_j-{\bomega}_j)_{S_{j*}}\|_2 \leq (1+\xi)s_0^{1/2} \|\widehat{\bomega}_j-{\bomega}_j\|_2.
\end{align*}
which together with \eqref{sbx} entails that
\begin{eqnarray}\label{sbxj}
\mathbb{P} \left\{\|\widehat{\bomega}_j-{\bomega}_j\|_1\leq Cs_0\lambda_j\right\}\geq 1- o(p^{1-\delta}).
\end{eqnarray}
Then by inequalities \eqref{sdd}, \eqref{sddr}, and \eqref{sbxj}, we finally have
\begin{eqnarray}\label{kongbu}
\mathbb{P} \left\{\|\bpsi_{\widehat{S}^c\backslash  \{j\}}^{(j)}(\widehat{\bomega}_j-{\bomega}_j)\|_2\leq C(s_0 n)^{1/2} \lambda_j \right\}\geq 1- o(p^{1-\delta}),
\end{eqnarray}
which entails that
\begin{eqnarray*}
\mathbb{P} \left\{n^{-1/2}\|\bpsi_{\widehat{S}^c\backslash  \{j\}}^{(j)}(\widehat{\bomega}_j-{\bomega}_j)\|_2\leq Cs_0^{1/2}\lambda_j\right\}\geq 1- o(p^{1-\delta}).
\end{eqnarray*}

We then derive the deviation bounds of $\|\widehat{\bomega}_j-{\bomega}_j\|_q$ for $q\in(1,2)$. Note that the two inequalities $\sum_{j=1}^p a_j\leq b_1$ and $\sum_{j=1}^p a_j^2\leq b_2$ for $a_j\geq 0$ will imply
$$\sum_{j=1}^p a_j^q=\sum_{j=1}^p a_j^{2-q}a_j^{2q-2}\leq\left(\sum_{j=1}^p a_j\right)^{2-q}\left( \sum_{j=1}^p a_j^2\right)^{q-1}\leq b_1^{2-q}b_2^{q-1}$$
for any $q\in (1,2)$. Thus, we have
$$\|\widehat{\bomega}_j-{\bomega}_j\|_q^q\leq \|\widehat{\bomega}_j-{\bomega}_j\|_1^{2-q}\cdot(\|\widehat{\bomega}_j-{\bomega}_j\|_2^2)^{q-1},$$
which together with \eqref{sbx} and \eqref{sbxj} gives
\begin{eqnarray*}
\mathbb{P} \left\{\|\widehat{\bomega}_j-{\bomega}_j\|_q\leq Cs_0^{1/q}\lambda_j\right\}\geq 1- o(p^{1-\delta})
\end{eqnarray*}
for any $q\in (1,2)$. It concludes the proof of Theorem \ref{theob}.

\subsection{Proof of Theorem \ref{theoa}}

Under the conditions of Theorem  \ref{theoa}, it is clear that the results in \emph{Part 1} and \emph{Part 2} of the proof of Proposition \ref{pro2} are still valid.
In what follows, we mainly present the proof when $j \in \widehat{S}^c$ and the same argument applies to $j \in \widehat{S}$ for a given set $\widehat{S}$. Throughout the proof,  the regularization parameters $\lambda_j=(1+\epsilon)\sqrt{(2\delta/(n\phi_{jj})) \log(p)}$ with $\epsilon=\widetilde{M}(\xi'+1)/\{\widetilde{M}^{\ast}(\xi'-1)\} - 1>0$, which are the same as those in the proof of Theorem \ref{theob}. In addition, $C$ is some positive constant which can vary from expression to expression.

\medskip

\noindent\emph{Part 1: Deviation bounds of $\|\widehat{\bpi}_{-j}^{(j)}-\bpi_{-j}^{(j)}\|_2$.}
In this part,  our analysis will be conditioning on the event $\mathcal{E}'\cap\mathcal{E}''$, where $\mathcal{E}'=\bigcap_{k \neq j}\{\widetilde{M}^{\ast}\leq \|\bx_k\|_2/\sqrt{n} \leq \widetilde{M}\}$  and $\mathcal{E}''=\bigcap_{k\in\widehat{S}^c\backslash  \{j\}}\{\widetilde{M}^{\ast}\leq \|\bpsi_k^{(j)}\|_2/\sqrt{n} \leq \widetilde{M}\}$, respectively.
For clarity, we start with some additional analyses.
Consider the conditional distribution
\begin{align}\label{yuyu0}
\bx_j=\bX_{-j} \bpi_{-j}^{(j)}+ \textbf{e}_j,
\end{align}
where $\textbf{e}_j=\bx_j- \mathbb{E}(\bx_j|\bX_{-j})$ is independent of $\bX_{-j}$ and follows the distribution $N(0,(1/\phi_{jj})\bI_n)$,  and the regression coefficient vector  $\bpi_{-j}^{(j)}=-\bPhi_{j,-j}^\top/\phi_{jj}$. Let  $S^j=\{1 \leq k \leq p: k\neq j, \phi_{jk}\neq 0\}$, $z_{\infty}^j=\|\bX_{-j}^\top(\bx_j-\bX_{-j}\bpi^{(j)}_{-j})/n\|_{\infty}$, and $J_j=\widehat{S}\cup S^j$.

Since $j \in \widehat{S}^c$, we can rewrite (\ref{yuyu0}) as
$$\bx_j=\bX_{\widehat{S}}{\bpi}^{(j)}_{\widehat{S}}+\bX_{\widehat{S}^c\backslash  \{j\}}{\bpi}^{(j)}_{\widehat{S}^c\backslash  \{j\}}+\textbf{e}_j.$$
In view of Proposition \ref{pro1}, we have 
\begin{align}\nonumber
&(\widehat{\bpi}^{(j)}_{\widehat{S}},\widehat{\bpi}^{(j)}_{\widehat{S}^c\backslash  \{j\}}) = \mathop{\arg\min}_{\bb_{\widehat{S}},\bb_{\widehat{S}^c\backslash  \{j\}}} \Big\{(2n)^{-1}\big\| \bx_j - \bX_{\widehat{S}} \bb_{\widehat{S}}-\bX_{\widehat{S}^c\backslash  \{j\}}\bb_{\widehat{S}^c\backslash  \{j\}}\big\|_2^2 + \lambda_j \sum_{k \in \widehat{S}^c\backslash  \{j\} } v_{jk} |b_k| \Big\},
\end{align}
whose KKT condition shows
\begin{align*}
\begin{cases}
\nonumber-\bX_{\widehat{S}}^\top(\bx_j-\bX_{\widehat{S}}\widehat{\bpi}^{(j)}_{\widehat{S}}-\bX_{\widehat{S}^c\backslash  \{j\}}\widehat{\bpi}^{(j)}_{\widehat{S}^c\backslash  \{j\}})=\bzero,\\  \nonumber
-\bX_{\widehat{S}^c\backslash  \{j\}}^\top(\bx_j-\bX_{\widehat{S}}\widehat{\bpi}^{(j)}_{\widehat{S}}-\bX_{\widehat{S}^c\backslash  \{j\}}\widehat{\bpi}^{(j)}_{\widehat{S}^c\backslash  \{j\}})=-n\lambda_j\widehat{\bkappa}^{(j)}_{\widehat{S}^c\backslash  \{j\}}.
\end{cases}
\end{align*}

Here the $k$th component of $\widehat{\bkappa}^{(j)}$ satisfies
$$\widehat{\kappa}_k^{(j)} \begin{cases}
= v_{jk} \cdot \text{sign}(\widehat{\pi}^{(j)}_k) & \text{if} \ \widehat{\pi}^{(j)}_k \ne 0, k \in \widehat{S}^c\backslash  \{j\}, \\
\in [-v_{jk},v_{jk}] &  \text{if} \ \widehat{\pi}^{(j)}_k = 0,  k \in \widehat{S}^c\backslash  \{j\},
\end{cases}
$$
where  $v_{jk} = \|\bpsi_k^{(j)}\|_2/\sqrt{n}$ and $\widehat{\pi}_k^{(j)}$ denotes the $k$th component of $\widehat{\bpi}^{(j)}$. 
Without loss of generality, we assume that $\widehat{S}=\{1,\ldots,d\}$. 
Thus, based on the above KKT condition, it follows that
\begin{eqnarray}\label{jiezheng}
\bX_{-j}^\top\bX_{-j}(\widehat{\bpi}_{-j}^{(j)}-\bpi_{-j}^{(j)})/n=\bX_{-j}^\top(\bx_j-\bX_{-j}\bpi^{(j)}_{-j})/n-\lambda_j\widehat{\bkappa}^{( j \ast)},
\end{eqnarray}
where $\widehat{\bkappa}^{( j \ast)}=(\bzero^\top,(\widehat{\bkappa}^{(j)}_{\widehat{S}^c\backslash  \{j\}})^\top)^\top$.

In what follows, we proceed to use similar arguments as those in \cite{ye2010} to construct bounds for $\|\widehat{\bpi}_{-j}^{(j)}-\bpi_{-j}^{(j)}\|_2$ and the main difference here is that the partially penalized regression problem can be more complex. To this end, we first analyze the properties of $\widehat{\bpi}_{-j}^{(j)}-\bpi_{-j}^{(j)}$.
Denote by $J_j^c = -j \backslash J_j = \{1,\ldots,p\}\backslash (J_j \cup \{j\})$. 
Thus, it follows from (\ref{jiezheng}) and $\bv^\top\bX_{-j}^\top(\bx_j-\bX_{-j}\bpi^{(j)}_{-j})/n\leq  z_{\infty}^j \|\bv\|_1$ that
\begin{align}\nonumber
\bv^\top\bX_{-j}^\top\bX_{-j}(\widehat{\bpi}_{-j}^{(j)}-\bpi_{-j}^{(j)})/n&=\bv^\top\bX_{-j}^\top(\bx_j-\bX_{-j}\bpi^{(j)}_{-j})/n-\lambda_j\bv^\top\widehat{\bkappa}^{( j \ast)}\\ \nonumber
&\leq z_{\infty}^j \|\bv\|_1-\lambda_j\bv^\top\widehat{\bkappa}^{( j \ast)}
\end{align}
for all vectors $\bv \in \mathbb{R}^{p-1}$. If $\bv$ satisfies $\text{sign}(\bv_{J_j^c})=\text{sign}((\widehat{\bpi}_{-j}^{(j)}-\bpi_{-j}^{(j)})_{J_j^c})=\text{sign}((\widehat{\bpi}_{-j}^{(j)})_{J_j^c})$, with the aid of $v_{jk}\geq \widetilde{M}^{\ast}$ for any $k \in \widehat{S}^c\backslash  \{j\}$, it follows that
\begin{align*}
-\bv^\top\widehat{\bkappa}^{( j \ast)}=-(\bv_{J_j}^\top\widehat{\bkappa}_{J_j}^{( j \ast)}+\bv_{J_j^c}^\top\widehat{\bkappa}_{J_j^c}^{( j \ast)})\leq -\bv_{J_j}^\top\widehat{\bkappa}_{J_j}^{( j \ast)}-\widetilde{M}^{\ast}\|\bv_{J_j^c}\|_1,
\end{align*}
which together with $\widehat{S}\subset J_j$ and  $v_{jk}\leq \widetilde{M}$ for any $k \in \widehat{S}^c\backslash  \{j\}$ further gives
\begin{align*}
-\bv^\top\widehat{\bkappa}^{( j \ast)}\leq \widetilde{M}\|\bv_{J_j}\|_1-\widetilde{M}^{\ast}\|\bv_{J_j^c}\|_1.
\end{align*}

Combining these results gives
\begin{align}\label{yuebing}
\bv^\top\bX_{-j}^\top\bX_{-j}(\widehat{\bpi}_{-j}^{(j)}-\bpi_{-j}^{(j)})/n\leq (z_{\infty}^j+ \widetilde{M}\lambda_j)\|\bv_{J_j}\|_1+(z_{\infty}^j- \widetilde{M}^{\ast}\lambda_j)\|\bv_{J_j^c}\|_1
\end{align}
for all vectors $\bv \in \mathbb{R}^{p-1}$ satisfying $\text{sign}(\bv_{J_j^c})=\text{sign}((\widehat{\bpi}_{-j}^{(j)}-\bpi_{-j}^{(j)})_{J_j^c})=\text{sign}((\widehat{\bpi}_{-j}^{(j)})_{J_j^c})$.
Based on the above inequality, if $z_{\infty}^j < \widetilde{M}^{\ast}\lambda_j$, it follows from $\bv=\widehat{\bpi}_{-j}^{(j)}-\bpi_{-j}^{(j)}$ that
$$(z_{\infty}^j+ \widetilde{M}\lambda_j)\|(\widehat{\bpi}_{-j}^{(j)}-\bpi_{-j}^{(j)})_{J_j}\|_1+(z_{\infty}^j- \widetilde{M}^{\ast}\lambda_j)\|(\widehat{\bpi}_{-j}^{(j)}-\bpi_{-j}^{(j)})_{J_j^c}\|_1\geq0,$$
which entails that
\begin{align}\label{zixidian}
\|(\widehat{\bpi}_{-j}^{(j)}-\bpi_{-j}^{(j)})_{J_j^c}\|_1 \leq \xi \|(\widehat{\bpi}_{-j}^{(j)}-\bpi_{-j}^{(j)})_{J_j}\|_1
\end{align}
for  $\xi\geq (\widetilde{M}\lambda_j+z_{\infty}^j)/(\widetilde{M}^{\ast}\lambda_j-z_{\infty}^j)$.

Moreover, in view of (\ref{jiezheng}), with $\textbf{e}^k= (0,\ldots,0,\underbrace{1}_{k\text{th}},0,\ldots,0)\in \mathbb{R}^{p-1}$, it follows from $\bv=(\widehat{\pi}^{(j)}_k -{\pi}^{(j)}_k )\textbf{e}^k$  that
\begin{align}\label{fanfan}
 (\widehat{\pi}^{(j)}_k -{\pi}^{(j)}_k )\bx_k^\top\bX_{-j}(\widehat{\bpi}_{-j}^{(j)}-\bpi_{-j}^{(j)})/n\leq z_{\infty}^j|\widehat{\pi}^{(j)}_k -{\pi}^{(j)}_k |-\widetilde{M}^{\ast}\lambda_j|\widehat{\pi}^{(j)}_k -{\pi}^{(j)}_k |
\end{align}
for any $ k\in J_j^c$, where the inequality can be derived by applying an argument similar to that in \eqref{yuebing}.
With the aid of \eqref{fanfan}, if  $z_{\infty}^j\leq\widetilde{M}^{\ast}\lambda_j$, it follows that
\begin{align}\label{henzixi}
 (\widehat{\pi}^{(j)}_k -{\pi}^{(j)}_k )\bx_k^\top\bX_{-j}(\widehat{\bpi}_{-j}^{(j)}-\bpi_{-j}^{(j)})/n \leq 0
\end{align}
for any $ k\in J_j^c$.
In view of \eqref{zixidian} and \eqref{henzixi}, if $z_{\infty}^j<\widetilde{M}^{\ast}\lambda_j$, we finally get
\begin{align}\label{chaozixi}
\widehat{\bpi}_{-j}^{(j)}-\bpi_{-j}^{(j)}\in \mathcal{G}_{-}(\xi,J_j)
 \end{align}
 for $\xi\geq (\widetilde{M}\lambda_j+z_{\infty}^j)/(\widetilde{M}^{\ast}\lambda_j-z_{\infty}^j)$, where the sign-restricted cone $\mathcal{G}_{-}(\xi,J_j)$  introduced in \cite{ye2010} is defined as
$$\mathcal{G}_{-}(\xi,J_j)=\{\bu\in \mathbb{R}^{p-1}:\|\bu_{J_j^c}\|_1 \leq \xi \|\bu_{J_j}\|_1, u_k\bx_k^\top\bX_{-j}\bu/n\leq 0, \forall k\in J_j^c\}.$$

In order to proceed, our analysis will be conditioning on the event $\mathcal{E}'''= \{z_{\infty}^j\leq \lambda_j\widetilde{M}^{\ast}(\xi'-1)/(\xi'+1)\}$ for some constant $\xi'>1$. Note that
$z_{\infty}^j\leq \lambda_j\widetilde{M}^{\ast}(\xi'-1)/(\xi'+1) < \widetilde{M}^{\ast}\lambda_j$
holds on this event, which together with \eqref{chaozixi} implies that
\begin{align*}
\widehat{\bpi}_{-j}^{(j)}-\bpi_{-j}^{(j)}\in \mathcal{G}_{-}(\xi,J_j)
 \end{align*}
for $\xi\geq (\widetilde{M}\lambda_j+z_{\infty}^j)/(\widetilde{M}^{\ast}\lambda_j-z_{\infty}^j)$. We then clarify that the magnitude of $(\widetilde{M}\lambda_j+z_{\infty}^j)/(\widetilde{M}^{\ast}\lambda_j-z_{\infty}^j)$ is around a constant level such that $\xi$ can be a constant.

In view of $z_{\infty}^j\leq \lambda_j\widetilde{M}^{\ast}(\xi'-1)/(\xi'+1)$, it follows that
$$\frac{\widetilde{M}\lambda_j+z_{\infty}^j}{\widetilde{M}^{\ast}\lambda_j-z_{\infty}^j}+1=\frac{\widetilde{M}\lambda_j+\widetilde{M}^{\ast}\lambda_j}{\widetilde{M}^{\ast}\lambda_j-z_{\infty}^j}\leq
\frac{\widetilde{M}\lambda_j+\widetilde{M}^{\ast}\lambda_j}{\widetilde{M}^{\ast}\lambda_j-\lambda_j\widetilde{M}^{\ast}(\xi'-1)/(\xi'+1)}=\frac{(\xi'+1)(\widetilde{M}+\widetilde{M}^{\ast})}{2\widetilde{M}^{\ast}},$$
which entails that
$$\frac{\widetilde{M}\lambda_j+z_{\infty}^j}{\widetilde{M}^{\ast}\lambda_j-z_{\infty}^j}\leq \frac{(\xi'+1)(\widetilde{M}+\widetilde{M}^{\ast})}{2\widetilde{M}^{\ast}}-1=1+\frac{\xi'-1}{2}+\frac{(\xi'+1)\widetilde{M}-2\widetilde{M}^{\ast}}{2\widetilde{M}^{\ast}},$$
which is a constant larger than one in view of $\xi' > 1$. Combining this inequality with $(\widetilde{M}\lambda_j+z_{\infty}^j)/(\widetilde{M}^{\ast}\lambda_j-z_{\infty}^j)\geq (\widetilde{M}\lambda_j)/(\widetilde{M}^{\ast}\lambda_j)\geq 1$ further gives
\begin{align}\label{yanzi}
1\leq \frac{\widetilde{M}\lambda_j+z_{\infty}^j}{\widetilde{M}^{\ast}\lambda_j-z_{\infty}^j}\leq 1+\frac{\xi'-1}{2}+\frac{(\xi'+1)\widetilde{M}-2\widetilde{M}^{\ast}}{2\widetilde{M}^{\ast}}.
\end{align}

Moreover, since $z_{\infty}^j\leq \lambda_j\widetilde{M}^{\ast}(\xi'-1)/(\xi'+1)\leq \lambda_j\widetilde{M}(\xi'-1)/(\xi'+1)$, it follows from (\ref{jiezheng}) that
\begin{align*}
\|\bX_{-j}^\top\bX_{-j}(\widehat{\bpi}_{-j}^{(j)}-\bpi_{-j}^{(j)})/n\|_{\infty}\leq \widetilde{M}\lambda_j+z_{\infty}^j
\leq \widetilde{M}\lambda_j+\lambda_j\widetilde{M}(\xi'-1)/(\xi'+1)=2\widetilde{M}\lambda_j\xi'/(\xi'+1).
\end{align*}
Thus, in view of $|J_j|\leq s_0+d$, combining the above results gives
\begin{align}\label{xuezhuan}
\|\widehat{\bpi}_{-j}^{(j)}-\bpi_{-j}^{(j)}\|_2\leq \frac{|J_j|^{1/2}(\widetilde{M}\lambda_j+z_{\infty}^j)}{F_2(\xi,J_j)}\leq\frac{\{2\widetilde{M}\xi'/(\xi'+1)\}(s_0+d)^{1/2}\lambda_j}{F_2(\xi,J_j)},
\end{align}
where the sign-restricted cone invertibility factor $F_2(\xi,S_j)$ introduced in \cite{ye2010} is defined as
\begin{align*}
  F_2 (\xi,J_j)  = \inf \left\{|J_j|^{1/2}\|\bX_{-j}^\top\bX_{-j} \bu/n\|_{\infty}/ \|\bu\|_{2}: \bu \in \mathcal{G}_{-}(\xi,J_j) \right\}. 
\end{align*}

Then we will derive the probability bound of the event $\mathcal{E}'''$, which can be obtained by making use of the inequality that
\begin{align}\label{wuyu}
\int_{x}^{\infty}e^{-t^2/2}dt \leq e^{-x^2/2}/x
\end{align}
for $x>0$. Since $\textbf{e}_j$ is independent of $\bX_{-j}$ and follows the distribution $N(0,(1/\phi_{jj})\bI_n)$, it yields that
$$\textbf{e}_j^\top{\bx}_{k}|{\bx}_{k}\sim N(0,\|\bx_k\|_2^2/\phi_{jj})$$
for any $k\neq j$, which together with \eqref{wuyu}  gives
$$\mathbb{P}\Big(\frac{\textbf{e}_j^\top{\bx}_{k}}{\sqrt{\|\bx_k\|_2^2/\phi_{jj}}}\big|{\bx}_{k} > \sqrt{2\delta \log(p)}\Big)\le  \frac{1}{\sqrt{2\pi}}\cdot \frac{p^{-\delta}}{\sqrt{2\delta \log(p)}}.$$

Note that the right hand side of the above inequality is independent of $\bx_k$. Thus, we can obtain
$$\mathbb{P}\Big(\frac{\textbf{e}_j^\top{\bx}_{k}}{\sqrt{\|\bx_k\|_2^2/\phi_{jj}}}> \sqrt{2\delta \log(p)}\Big)\le  \frac{1}{\sqrt{2\pi}}\cdot \frac{p^{-\delta}}{\sqrt{2\delta \log(p)}}.$$
Based on this inequality and taking the probability of the event $\mathcal{E}'$ into consideration, with $\lambda_j=[\widetilde{M}(\xi'+1)/\{\widetilde{M}^{\ast}(\xi'-1)\}] \sqrt{(2\delta/(n\phi_{jj})) \log(p)}$, it follows that
\begin{align}\nonumber
&1-\mathbb{P}\left(z_{\infty}^j\leq \frac{\lambda_j\widetilde{M}^{\ast}(\xi'-1)}{\xi'+1}\right)
=\mathbb{P}\Big(\max_{k\neq j, 1\le k\le p}\frac{|\textbf{e}_j^\top{\bx}_{k}|}{n}> \widetilde{M}\sqrt{\frac{2\delta\log(p)}{n\phi_{jj}} }\Big) \\\nonumber
\le \ & \mathbb{P}\Big(\max_{k\neq j, 1\le k\le p}\frac{|\textbf{e}_j^\top{\bx}_{k}|}{n}> \widetilde{M}\sqrt{\frac{2\delta\log(p)}{n\phi_{jj}} }\Big| \mathcal{E}'\Big)+\mathbb{P}((\mathcal{E}')^c) \\\nonumber
\le \ &  2 \sum_{k\neq j, 1\le k\le p} \mathbb{P}\Big(\frac{\textbf{e}_j^\top{\bx}_{k}}{\sqrt{\|\bx_k\|_2^2/\phi_{jj}}} > \frac{\widetilde{M}\sqrt{n}}{\|\bx_k\|_2}\sqrt{2\delta \log(p)}\Big| \mathcal{E}'\Big) +o(p^{1-\delta})\\\label{luy1}
\leq \ & 2 \sum_{k\neq j, 1\le k\le p} \mathbb{P}\Big(\frac{\textbf{e}_j^\top{\bx}_{k}}{\sqrt{\|\bx_k\|_2^2/\phi_{jj}}} > \sqrt{2\delta \log(p)}\Big) +o(p^{1-\delta})
\le  o(p^{1-\delta}).
\end{align}

We still need to know the magnitude of $F_2(\xi,J_j)$. Since the eigenvalues of $\bSigma$ are bounded within the interval $[1/L, L]$ for some constant $L\geq 1$, we can deduce that the eigenvalues of $\bSigma_{-j,-j}$ are also bounded within the interval $[1/L, L]$ by Lemma \ref{sigma}. Since $s_0+d=o(n/\log (p))$ can be derived by Condition \ref{weidu} and $d=o(n/\log (p))$, in view of Lemma \ref{lem2}, there are positive constants $c_1$, $c_2$, $c_3$ such that
\begin{align*}
F_2 (\xi,J_j)  \geq c_1
\end{align*}
with probability at least $1- c_2\exp(-c_3n)$. Since $\log(p)/n=o(1)$, it follows that
$$\frac{\exp(-c_3n)}{p^{1-\delta}}=\frac{p^{\delta-1}}{\exp(c_3n)}=\exp\{(\delta-1)\log p -c_3n\}=o(1),$$
which entails $c_2\exp(-c_3n)=o(p^{1-\delta})$.

Moreover, with the aid of (\ref{ycj}), (\ref{ycj1}), and \eqref{luy1}, we know that
\begin{align*}
&\mathbb{P}(\mathcal{E}'\cap\mathcal{E}''\cap\mathcal{E}''')=1-\mathbb{P}\left((\mathcal{E}')^c\cup(\mathcal{E}'')^c\cup(\mathcal{E}''')^c\right)\\
\geq & 1-\mathbb{P}((\mathcal{E}')^c)-\mathbb{P}((\mathcal{E}'')^c)-\mathbb{P}((\mathcal{E}''')^c)=1-o(p^{1-\delta}).
\end{align*}
Thus, combining these results with \eqref{xuezhuan} gives
\begin{eqnarray}\label{kouzhao1}
\mathbb{P} \left\{\|\widehat{\bpi}_{-j}^{(j)}-\bpi_{-j}^{(j)}\|_2\leq C(s_0+d)^{1/2}\lambda_j\right\}\geq 1-o(p^{1-\delta}).
\end{eqnarray}

\medskip
\noindent \emph{Part 2: Deviation bounds of $n^{-1/2}\|\bX_{-j} (\widehat{\bpi}_{-j}^{(j)}-\bpi_{-j}^{(j)})\|_2$ and $\|\widehat{\bpi}_{-j}^{(j)}-\bpi_{-j}^{(j)}\|_q$ with $q\in[1,2)$.}
We first construct bounds for $\|\bX_{-j} (\widehat{\bpi}_{-j}^{(j)}-\bpi_{-j}^{(j)})\|_2$. Some simple algebra shows that
\begin{align}\label{sbd}
\|\bX_{-j} (\widehat{\bpi}_{-j}^{(j)}-\bpi_{-j}^{(j)})\|_2^2\leq \|\bX_{-j}^\top\bX_{-j}(\widehat{\bpi}_{-j}^{(j)}-\bpi_{-j}^{(j)})\|_{\infty}\|\widehat{\bpi}_{-j}^{(j)}-\bpi_{-j}^{(j)}\|_1.
\end{align}
In the event $\mathcal{E}'''=\{z_{\infty}^j\leq \lambda_j\widetilde{M}^{\ast}(\xi'-1)/(\xi'+1)\}$, it follows that
\begin{align*}
\|\bX_{-j}^\top\bX_{-j}(\widehat{\bpi}_{-j}^{(j)}-\bpi_{-j}^{(j)})/n\|_{\infty}\leq \widetilde{M}\lambda_j+z_{\infty}^j \leq  \widetilde{M}\lambda_j+\lambda_j\widetilde{M}^{\ast}(\xi'-1)/(\xi'+1),
\end{align*}
which further yields
\begin{align}\label{sbdr}\|\bX_{-j}^\top\bX_{-j}(\widehat{\bpi}_{-j}^{(j)}-\bpi_{-j}^{(j)})\|_{\infty}\leq n \{\widetilde{M}+\widetilde{M}^{\ast}(\xi'-1)/(\xi'+1)\}\lambda_j.
\end{align}

Moreover, since $\widehat{\bpi}_{-j}^{(j)}-\bpi_{-j}^{(j)}\in \mathcal{G}_{-}(\xi,J_j)$, we have
\begin{align*}
\|\widehat{\bpi}_{-j}^{(j)}-\bpi_{-j}^{(j)}\|_1&=\|(\widehat{\bpi}_{-j}^{(j)}-\bpi_{-j}^{(j)})_{J_j^c}\|_1 + \|(\widehat{\bpi}_{-j}^{(j)}-\bpi_{-j}^{(j)})_{J_j}\|_1\leq (1+\xi) \|(\widehat{\bpi}_{-j}^{(j)}-\bpi_{-j}^{(j)})_{J_j}\|_1\\
&\leq (1+\xi)(s_0+d)^{1/2} \|(\widehat{\bpi}_{-j}^{(j)}-\bpi_{-j}^{(j)})_{J_j}\|_2 \leq (1+\xi)(s_0+d)^{1/2} \|\widehat{\bpi}_{-j}^{(j)}-\bpi_{-j}^{(j)}\|_2,
\end{align*}
which together with \eqref{kouzhao1} entails that
\begin{eqnarray}\label{kouzhao2}
\mathbb{P} \left\{\|\widehat{\bpi}_{-j}^{(j)}-\bpi_{-j}^{(j)}\|_1\leq C(s_0+d)\lambda_j\right\}\geq 1- o(p^{1-\delta}).
\end{eqnarray}
Then with the aid of \eqref{sbd}, \eqref{sbdr}, and \eqref{kouzhao2}, we finally get
\begin{eqnarray}\label{laoxu2}
\mathbb{P} \left\{\|\bX_{-j} (\widehat{\bpi}_{-j}^{(j)}-\bpi_{-j}^{(j)})\|_2\leq C\{(s_0+d) n\}^{1/2} \lambda_j \right\}\geq 1- o(p^{1-\delta}),
\end{eqnarray}
which further entails that
\begin{eqnarray*}
\mathbb{P} \left\{n^{-1/2}\|\bX_{-j} (\widehat{\bpi}_{-j}^{(j)}-\bpi_{-j}^{(j)})\|_2\leq C(s_0+d)^{1/2}\lambda_j\right\}\geq 1- o(p^{1-\delta}).
\end{eqnarray*}

We can also derive the deviation bounds of $\|\widehat{\bpi}_{-j}^{(j)}-\bpi_{-j}^{(j)}\|_q$ with $q\in(1,2)$. Note that the two inequalities $\sum_{j=1}^p a_j\leq b_1$ and $\sum_{j=1}^p a_j^2\leq b_2$ for $a_j\geq 0$ will imply
$$\sum_{j=1}^p a_j^q=\sum_{j=1}^p a_j^{2-q}a_j^{2q-2}\leq\left(\sum_{j=1}^p a_j\right)^{2-q}\left( \sum_{j=1}^p a_j^2\right)^{q-1}\leq b_1^{2-q}b_2^{q-1}$$
for any $q\in (1,2)$. Thus, we have
$$\|\widehat{\bpi}_{-j}^{(j)}-\bpi_{-j}^{(j)}\|_q^q\leq \|\widehat{\bpi}_{-j}^{(j)}-\bpi_{-j}^{(j)}\|_1^{2-q}\cdot(\|\widehat{\bpi}_{-j}^{(j)}-\bpi_{-j}^{(j)}\|_2^2)^{q-1},$$
which together with \eqref{kouzhao1} and \eqref{kouzhao2} gives
\begin{eqnarray*}
\mathbb{P} \left\{\|\widehat{\bpi}_{-j}^{(j)}-\bpi_{-j}^{(j)}\|_q\leq C(s_0+d)^{1/q}\lambda_j\right\}\geq 1- o(p^{1-\delta})
\end{eqnarray*}
for any $q\in (1,2)$. It completes the proof of Theorem \ref{theoa}.

\subsection{Proof of Theorem \ref{theo1}}

Under the conditions of Theorem  \ref{theo1}, it is clear that the results in \emph{Part 1} and \emph{Part 2} of the proof of Proposition \ref{pro2} are still valid. Denote by $\mathcal{E}=\big\{S_1\subset\widehat{S}\big\}$ the event in Definition \ref{C2}. By Definition \ref{C2}, $\widehat{S}$ is independent of the data $(\bX, \by)$, so we treat $\widehat{S}$ as given. Conditional on the event $\mathcal{E}$, we will first analyze the properties of some key statistics before deriving the confidence intervals of $\beta_j$. In what follows, we mainly present the proof when $j \in \widehat{S}^c$ and the same argument applies to $j \in \widehat{S}$. Throughout the proof, $C$ is some positive constant which can vary from expression to expression.


\medskip
\noindent \emph{Part 1: Deviation bounds of $\|\bz_j\|_2$}. We derive the deviation bounds of $\|\bz_j\|_2$ in two different cases that correspond to the asymptotic results in Theorems \ref{theob} and \ref{theoa}, respectively. In each case, the assumption on $s_1$ is different but we get asymptotically equivalent deviation bounds of $\|\bz_j\|_2$. In both cases, the regularization parameters $\lambda_j=(1+\epsilon)\sqrt{(2\delta/(n\phi_{jj})) \log(p)}$ with $\epsilon=\widetilde{M}(\xi'+1)/\{\widetilde{M}^{\ast}(\xi'-1)\} - 1>0$, which are the same as those in the proof of Theorem \ref{theoa}.

\emph{Case 1: $s_1 = o(n/s_0)$.} Since $s_0 = o(n/\log (p))$ and $d=O(s_1) = o(n/s_0)$, the convergence rates established in Theorem \ref{theob} are asymptotically vanishing. Then in this case, we mainly utilize the results in Theorem \ref{theob}.
Since $\textbf{e}_j$ is independent of $\bX_{\widehat{S}}$, it follows from $\textbf{e}_j^*=(\bI_n-\mathbb{\bP}_{\widehat{S}})\textbf{e}_j$ defined in \eqref{lud} that
$$\|\textbf{e}_j^*\|_2^2|\bX_{\widehat{S}}\sim \chi_{(n-d)}^2.$$
In view of $d =o(n)$, applying the same argument as that in (\ref{changdu}) gives
\begin{eqnarray*}
\mathbb{P}\left\{\widetilde{M}^{\ast}_4\leq (\|\textbf{e}_j^*\|_2/\sqrt{n}) \leq \widetilde{M}_4 \ | \ \bX_{\widehat{S}} \right\}\geq 1- o(p^{1-\delta}),
\end{eqnarray*}
where $\widetilde{M}_4\geq \widetilde{M}^{\ast}_4$  are some positive constants. Since the right hand side of the above inequality is independent of $\bX_{\widehat{S}}$, we further have
\begin{eqnarray}\label{duziteng}
\mathbb{P}\left\{\widetilde{M}^{\ast}_4\leq \|\textbf{e}_j^*\|_2/\sqrt{n}\leq \widetilde{M}_4\right\}\geq 1- o(p^{1-\delta}).
\end{eqnarray}

In addition, it follows from $\bpsi_j^{(j)}=\bpsi_{\widehat{S}^c\backslash  \{j\}}^{(j)} \bomega_j+\textbf{e}_j^*$ and $\bz_j = \bpsi_j^{(j)} - \bpsi_{\widehat{S}^c\backslash  \{j\}}^{(j)} \widehat{\bomega}_j$ that
$$\textbf{e}_j^*-\bz_j=\bpsi_{\widehat{S}^c\backslash  \{j\}}^{(j)}(\widehat{\bomega}_j-{\bomega}_j),$$
which together with inequality (\ref{kongbu}) gives
\begin{equation}\label{wuyul}
\mathbb{P}\left\{\|\textbf{e}_j^*\|_2-C(s_0n)^{1/2}\lambda_j \leq\|\bz_j\|_2\leq \|\textbf{e}_j^*\|_2+C(s_0n)^{1/2}\lambda_j \right\}\geq 1-  o(p^{1-\delta}).
\end{equation}
In view of $s_0 (\log p)/n = o(1)$, the above inequality together with the inequality (\ref{duziteng}) yields that
\begin{eqnarray}\label{tianzi1}
\mathbb{P}\left\{\|\bz_j\|_2\asymp n^{1/2}\right\}\geq 1-  o(p^{1-\delta}).
\end{eqnarray}
Therefore, in both cases, we have $\|\bz_j\|_2\asymp n^{1/2}$ with significant probability at least $1 - o(p^{1-\delta})$.

\emph{Case 2: $s_1 = o(n/\log (p))$.} In view of $d=O(s_1)$ by Definition \ref{C2}, it follows from $s_1 \log (p)/n = o(1)$ that $d \log (p)/n = o(1)$, so that the convergence rates in Theorem \ref{theoa} are asymptotically vanishing. 
Then in this case, we mainly utilize the asymptotic results in Theorem \ref{theoa}. 
For $\textbf{e}_j=\bx_j-\bX_{-j} \bpi_{-j}^{(j)}$ defined in (\ref{yuyu0}), applying the same argument as that in (\ref{changdu}) gives
\begin{eqnarray}\label{laoxu3}
\mathbb{P}\left\{\widetilde{M}_3^{\ast}\leq\|\textbf{e}_j\|_2/\sqrt{n} \leq \widetilde{M}_3\right\}\geq 1- o(p^{1-\delta}),
\end{eqnarray}
where $\widetilde{M}_3\geq \widetilde{M}^{\ast}_3$  are some positive constants.
By  Proposition \ref{pro1}, we have $\bz_j=\bx_j-\bX_{-j}\widehat{\bpi}_{-j}^{(j)}$. It yields that 
$$\textbf{e}_j-\bz_j=\bX_{-j} (\widehat{\bpi}_{-j}^{(j)}-\bpi_{-j}^{(j)}),$$
which together with the inequality (\ref{laoxu2}) gives
\begin{equation}\label{zee}
\mathbb{P}\left\{\|\textbf{e}_j\|_2-\widetilde{M}\{(s_0+d)\log p\}^{1/2} \leq\|\bz_j\|_2\leq \|\textbf{e}_j\|_2+\widetilde{M}\{(s_0+d)\log p\}^{1/2} \right\}\geq 1-  o(p^{1-\delta}).
\end{equation}
In view of $s_0 (\log p)/n = o(1)$ and $d (\log p)/n = o(1)$, the above inequality together with inequality (\ref{laoxu3}) shows
\begin{eqnarray}\label{laowu4}
\mathbb{P}\left\{\|\bz_j\|_2\asymp n^{1/2}\right\}\geq 1-  o(p^{1-\delta}).
\end{eqnarray}

\medskip
\noindent \emph{Part 2: Deviation bounds of $\tau_j$}.
First of all, it follows from the equality \eqref{zxb} that
\begin{equation}\label{KKT}
\begin{aligned}
\bz_j^\top \bx_j &=\bz_j^\top \bpsi_j^{(j)}=\|\bz_j\|_2^2 + (\bpsi_{\widehat{S}^c\backslash  \{j\}}^{(j)} \widehat{\bomega}_j)^\top \bz_j = \|\bz_j\|_2^2 + \sqrt{n}\lambda_j \sum_{k \in \widehat{S}^c\backslash  \{j\}} (\|\bpsi_k^{(j)}\|_2 \cdot |\homega_{j,k}|) \geq \|\bz_j\|_2^2,
\end{aligned}
\end{equation}
where the third equality above follows from the Karush-Kuhn-Tucker (KKT) condition of the Lasso estimator, which gives
\begin{equation}\label{KKT'}
\bz_j^\top \bpsi_k^{(j)} =(\bpsi_j^{(j)} - \bpsi_{\widehat{S}^c\backslash  \{j\}}\widehat{\bomega}_j)^\top\bpsi_k^{(j)} = \sqrt{n}\lambda_j \|\bpsi_k^{(j)}\|_2 \cdot \sgn(\homega_{j,k})
\end{equation}
with $\homega_{j,k}$ being the $k$th component of $\widehat{\bomega}_{j}$, for any $k \in A = \{k \in \widehat{S}^c\backslash  \{j\}: \sgn(\homega_{j,k}) \neq 0\}$.
By (\ref{KKT}), we immediately have
\begin{equation} \label{taud}
\begin{aligned}
\tau_j &=\|\bz_j\|_2/|\bz_j^\top \bx_j| \leq 1/\|\bz_j\|_2.
\end{aligned}
\end{equation}

It remains to find the lower bound of $\tau_j$. With the aid of (\ref{ycj1}), it follows that with probability at least $1 - o(p^{1-\delta})$,
\begin{equation}\label{KKT1}
\begin{aligned}
\bz_j^\top \bx_j &=\bz_j^\top \bpsi_j^{(j)} = \|\bz_j\|_2^2 + \sqrt{n}\lambda_j \sum_{k \in \widehat{S}^c\backslash  \{j\}} (\|\bpsi_k^{(j)}\|_2 \cdot |\homega_{j,k}|) \leq\|\bz_j\|_2^2+n \widetilde{M}_2\lambda_j\|\widehat{\bomega}_j\|_1.
\end{aligned}
\end{equation}
Base on (\ref{taud}) and (\ref{KKT1}), it follows that
\begin{eqnarray}\label{yqwl}
\mathbb{P}\left\{\frac{ \|\bz_j\|_2}{\|\bz_j\|_2^2+ \widetilde{M}_2n \lambda_j\|\widehat{\bomega}_j\|_1}\leq\tau_j\leq \frac{1}{\|\bz_j\|_2 }\right\}\geq 1- o(p^{1-\delta}).
\end{eqnarray}

We then derive the bounds of $\|\widehat{\bomega}_j\|_1$. Note that by the triangle inequality,
$$\|\widehat{\bomega}_j\|_1\leq \|\bomega_j\|_1+\|\widehat{\bomega}_j-\bomega_j\|_1.$$
According to the equality \eqref{lud}, $\bomega_j={\bpi}^{(j)}_{\widehat{S}^c\backslash  \{j\}}=\bPhi_{j,{\widehat{S}^c\backslash  \{j\}}}^\top/\phi_{jj}$. Thus, we have $\|\bomega_j\|_1\leq \|\bpi_{-j}^{(j)}\|_1\leq s_0^{1/2}\|\bpi_{-j}^{(j)}\|_2$, where $\bpi_{-j}^{(j)}=-\bPhi_{j,-j}^\top/\phi_{jj}$. 
Since the eigenvalues of $\bPhi$ are bounded within the interval $[1/L, L]$ for some constant $L\geq 1$, it follows from Lemma \ref{sigma} that
$$1/L\leq \|\bPhi_{j,-j}\|_2\leq L,  \ \  1/L\leq \sigma_{jj}\leq L,  \ \ \text{and} \ \ \sigma_{jj}\phi_{jj}\geq 1.$$
Based on these facts, it holds that
 $$\|\bpi_{-j}^{(j)}\|_2\leq \|\bPhi_{j,-j}\|_2/\phi_{jj}\leq L\sigma_{jj}\leq L^2=O(1),$$
which further gives
$$\|\bomega_j\|_1\leq \|\bpi_{-j}^{(j)}\|_1\leq  s_0^{1/2}\|\bpi_{-j}^{(j)}\|_2\leq O(s_0^{1/2}).$$

Since $s_0(\log p)/n= o(1)$, combining these results with \eqref{sbxj} gives
\begin{eqnarray*}
\mathbb{P}\left\{\|\widehat{\bomega}_j\|_1\leq \|\bomega_j\|_1+\|\widehat{\bomega}_j-\bomega_j\|_1=O( s_0^{1/2})+O(s_0\lambda_j)=O( s_0^{1/2})\right\}\geq 1- o(p^{1-\delta}).
\end{eqnarray*}
Utilizing $s_0(\log p)/n= o(1)$ again, we get $n \lambda_j\|\widehat{\bomega}_j\|_1 = o(n)$ based on the bound of $\|\widehat{\bomega}_j\|_1$ in the above inequality. This along with (\ref{tianzi1}) and \eqref{yqwl} gives that
\begin{eqnarray*}
\mathbb{P}\left\{\tau_j\asymp n^{-1/2}\right\}\geq 1- o(p^{1-\delta}).
\end{eqnarray*}

Moreover, it is not difficult to see that with probability at least $1 - o(p^{1-\delta})$,
\[\lim_{n \to \infty} \tau_j n^{1/2} = \lim_{n \to \infty} n^{1/2} \|\bz_j\|_2/|\bz_j^\top \bx_j| = \lim_{n \to \infty} n^{1/2} / \|\bz_j\|_2,\]
where the last equality is due to the previously established inequalities (\ref{KKT}) and (\ref{KKT1}) such that
$$\|\bz_j\|_2^2 \leq |\bz_j^\top \bx_j|\leq \|\bz_j\|_2^2+ Cn \lambda_j\|\widehat{\bomega}_j\|_1 = \|\bz_j\|_2^2 + o(n).$$
We will then derive the limit of $\tau_j n^{1/2}$ when $n \to \infty$ by showing that
$$\lim_{n \to \infty} n^{1/2} / \|\bz_j\|_2= \lim_{n \to \infty} n^{1/2} / \|\textbf{e}_j\|_2.$$

Since the above limit already holds in view of \eqref{zee} in \emph{Case 1} of \emph{Part 1}, we mainly prove it for \emph{Case 2} of \emph{Part 1}, where we have
\begin{eqnarray*}
\mathbb{P}\left\{\|\textbf{e}_j^*\|_2-C\{s_0\log p\}^{1/2} \leq\|\bz_j\|_2\leq \|\textbf{e}_j^*\|_2+C\{s_0\log p\}^{1/2} \right\}\geq 1-  o(p^{1-\delta})
\end{eqnarray*}
according to \eqref{wuyul}. We first derive the deviation of $\|\textbf{e}_j^*\|_2$ from  $\|\textbf{e}_j\|_2$.
Note that $\textbf{e}_j^*=(\bI_n-\mathbb{\bP}_{\widehat{S}})\textbf{e}_j$ is defined in \eqref{lud}. Since $\textbf{e}_j$ is independent of $\bX_{\widehat{S}}$,  it follows that
$$\|\mathbb{\bP}_{\widehat{S}}\textbf{e}_j\|_2^2|\bX_{\widehat{S}}\sim \chi_{(d)}^2.$$
Thus, applying the same argument as that in (\ref{changdu}) gives
\begin{eqnarray*}
\mathbb{P}\left\{\widetilde{M}^{\ast}_5\leq (\|\bP_{\widehat{S}}\textbf{e}_j\|_2/\sqrt{d})|\bX_{\widehat{S}} \leq \widetilde{M}_5\right\}\geq 1- o(p^{1-\delta}),
\end{eqnarray*}
where $\widetilde{M}_5\geq \widetilde{M}^{\ast}_5$ are two positive constants. Since the right hand side of the above inequality is independent of $\bX_{\widehat{S}}$, we get
\begin{eqnarray*}
\mathbb{P}\left\{\widetilde{M}^{\ast}_5\leq \|\bP_{\widehat{S}}\textbf{e}_j\|_2/\sqrt{d}\leq \widetilde{M}_5\right\}\geq 1- o(p^{1-\delta}).
\end{eqnarray*}

This together with $\|\textbf{e}_j\|_2-\|\mathbb{\bP}_{\widehat{S}}\textbf{e}_j\|_2\leq\|(\bI_n-\mathbb{\bP}_{\widehat{S}})\textbf{e}_j\|_2 = \|\textbf{e}_j^*\|_2 \leq \|\textbf{e}_j\|_2+\|\mathbb{\bP}_{\widehat{S}}\textbf{e}_j\|_2$ further gives
$$
\mathbb{P}\left\{\|\textbf{e}_j\|_2-\widetilde{M}^{\ast}_5\sqrt{d} \leq \|\textbf{e}_j^*\|_2\leq \|\textbf{e}_j\|_2+\widetilde{M}_5\sqrt{d} \right\}\geq 1-  o(p^{1-\delta}).
$$
Combining this result with \eqref{wuyul} yields that
\begin{eqnarray*}
\|\textbf{e}_j\|_2-\widetilde{M}^{\ast}_5\sqrt{d}-C\{s_0\log p\}^{1/2} \leq\|\bz_j\|_2\leq  \|\textbf{e}_j\|_2+\widetilde{M}_5\sqrt{d}+C\{s_0\log p\}^{1/2}
\end{eqnarray*}
holds with probability at least $1 - o(p^{1-\delta})$.

Since $d=o(n)$ and $s_0 (\log p)/n = o(1)$, we finally have
\[\lim_{n \to \infty} n^{1/2} / \|\bz_j\|_2= \lim_{n \to \infty} n^{1/2} / \|\textbf{e}_j\|_2 = \phi_{jj}^{-1/2},\]
where the last equality holds due to the facts that $\textbf{e}_j=\bx_j- \mathbb{E}(\bx_j|\bX_{-j})$ and $\textbf{e}_j$ has the distribution $N(\bzero,(1/\phi_{jj})\bI_n)$. It follows that
\[\lim_{n \to \infty} \tau_j n^{1/2}  = \lim_{n \to \infty} n^{1/2} / \|\bz_j\|_2 = \phi_{jj}^{-1/2}\]
holds with probability at least $1 - o(p^{1-\delta})$.

\medskip

\noindent \emph{Part 3: Deviation bounds of the inference bias $\Delta_j$.} Note that the inference bias $\Delta_j$ here is defined as
\begin{eqnarray*}
\Delta_j=\sum_{k \in \widehat{S}^c\backslash  \{j\}}\frac{\bz_j^\top\bx_k\beta_k}{\bz_j^\top\bx_j}.
\end{eqnarray*}
Then we have
\begin{eqnarray*}
|\Delta_j|=\Big|\sum_{k \in \widehat{S}^c\backslash  \{j\}}\frac{\bz_j^\top\bx_k\beta_k}{\bz_j^\top\bx_j}\Big|\leq \left(\max_{k \in \widehat{S}^c\backslash  \{j\}}\Big|\frac{\bz_j^\top\bx_k}{\bz_j^\top\bx_j}\Big|\right)\sum_{k \in {\widehat{S}^c}}|\beta_k|.
\end{eqnarray*}

Similar to \eqref{KKT} and \eqref{KKT'}, by the KKT condition and the equality \eqref{zxb}, for any $k\in \widehat{S}^c\backslash  \{j\}$, we can get
\[|\bz_j^\top\bx_k|=|\bz_j^\top\bpsi_k^{(j)}| = | (\bpsi_j^{(j)} - \bpsi_{\widehat{S}^c\backslash  \{j\}}\widehat{\bomega}_j)^\top\bpsi_k^{(j)}| \leq \sqrt{n}\lambda_j \|\bpsi_k^{(j)}\|_2.\]
Recall that the inequality (\ref{changdu1}) gives
\begin{equation*}
\begin{aligned}
&\mathbb{P}(\max_{k \in \widehat{S}^c\backslash  \{j\}} \frac{\|\bpsi_k^{(j)}\|_2}{\sqrt{n}} > \widetilde{M}_2)\leq \sum_{k \in \widehat{S}^c\backslash  \{j\}} \mathbb{P}(\frac{\|\bpsi_k^{(j)}\|_2}{\sqrt{n}} > \widetilde{M}_2)\leq p \cdot 4p^{-2\delta} = 4p^{1 - 2\delta} = o(p^{1 - \delta}).
\end{aligned}
\end{equation*}

Thus, with the aid of (\ref{laowu4}) and (\ref{KKT}), we have
\begin{eqnarray*}
\mathbb{P}\left\{\max_{k \in \widehat{S}^c\backslash  \{j\}}\Big|\frac{\bz_j^\top\bx_k}{\bz_j^\top\bx_j}\Big|\leq \max_{k \in \widehat{S}^c\backslash  \{j\}}\frac{\sqrt{n}\lambda_j \|\bpsi_k^{(j)}\|_2}{\|\bz_j\|_2^2} \leq C\sqrt{\log(p)/n}\right\}\geq 1- o(p^{1-\delta}).
\end{eqnarray*}
Since $S_1\subset\widehat{S}$ on the event $\mathcal{E}$, the above inequality together with $\|\bbeta_{S_1^c}\|_1=o(1/\sqrt{\log p})$ in Condition \ref{xinhao} entails that
\begin{eqnarray*}
\mathbb{P}\left\{|\Delta_j|\leq C\sqrt{\log(p)/n} \cdot o(1/\sqrt{\log p})  =o(n^{-1/2}) \right\}\geq 1- o(p^{1-\delta}).
\end{eqnarray*}

\medskip

\noindent \emph{Part 4: Confidence intervals of the coefficients $\beta_j$.}

The proof of this part is the same as that of Theorem \ref{haolei}, and, therefore, has been omitted.

\section{Additional technical details}\label{EC2} 

This part presents the proofs of the lemmas that are used in the proofs of the main theoretical results.



\subsection{Proof of Lemma \ref{tongfb}}

We first consider the conditional distribution
\begin{align*}
& \bx_j=\bX_{-j} \bpi_{-j}^{(j)}+ \textbf{e}_j,  
\end{align*}
where $\textbf{e}_j=\bx_j- \mathbb{E}(\bx_j|\bX_{-j})\sim N(\bzero,(1/\phi_{jj})\bI_n)$ and $\bpi_{-j}^{(j)}=-\bPhi_{j,-j}^\top/\phi_{jj}$.
Without loss of generality, we write $\bSigma$ as
$$\left(\begin{array}{cc}{\bSigma_{\widehat{S}\widehat{S}}}   & {\bSigma_{\widehat{S}\widehat{S}^c}} \\   {\bSigma_{\widehat{S}^c\widehat{S}}} & {\bSigma_{\widehat{S}^c\widehat{S}^c}}\end{array}\right).$$
By some simple algebra, we can get
\begin{align}\label{jige}
\bPhi = \bSigma^{-1}=\left(\begin{array}{cc}{\bSigma_{\widehat{S}\widehat{S}}}   & {\bSigma_{\widehat{S}\widehat{S}^c}} \\   {\bSigma_{\widehat{S}^c\widehat{S}}} & {\bSigma_{\widehat{S}^c\widehat{S}^c}}\end{array}\right)^{-1}=\left(\begin{array}{cc}{\ast}   & {\ast} \\   {\ast} & {\bSigma_{\widehat{S}^c\widehat{S}^c.1}^{-1}}\end{array}\right),
\end{align}
where $\bSigma_{\widehat{S}^c\widehat{S}^c.1}=\bSigma_{\widehat{S}^c\widehat{S}^c}-\bSigma_{\widehat{S}^c\widehat{S}}\bSigma_{\widehat{S}\widehat{S}}^{-1}\bSigma_{\widehat{S}\widehat{S}^c}$.

On the other hand, based on the conditional distribution of $\bX_{\widehat{S}^c}$ given $\bX_{\widehat{S}}$, it yields that
$$\brho^{(j)}_{\widehat{S}^c}\sim N(\bzero, \bI_n \otimes (\bSigma_{\widehat{S}^c\widehat{S}^c}-\bSigma_{\widehat{S}^c\widehat{S}}\bSigma_{\widehat{S}\widehat{S}}^{-1}\bSigma_{\widehat{S}\widehat{S}^c})).$$
Then it follows from \eqref{jige} that $$(\bSigma_{\widehat{S}^c\widehat{S}^c}-\bSigma_{\widehat{S}^c\widehat{S}}\bSigma_{\widehat{S}\widehat{S}}^{-1}\bSigma_{\widehat{S}\widehat{S}^c})^{-1}=\bSigma^{-1}_{\widehat{S}^c\widehat{S}^c.1}=\bPhi_{\widehat{S}^c\widehat{S}^c}.$$

Therefore, we immediately have
$$\brho^{(j)}_{\widehat{S}^c}\sim N(\bzero, \bI_n \otimes \bPhi_{\widehat{S}^c\widehat{S}^c}^{-1}),$$
which further yields
$$\brho^{(j)}_{\widehat{S}^c\backslash  \{j\}}\sim N(\bzero, \bI_n \otimes \bPhi_{\widehat{S}^c\backslash  \{j\}\widehat{S}^c\backslash  \{j\}}^{-1}).$$
It completes the proof of Lemma \ref{tongfb}.

\subsection{Proof of Lemma \ref{sigma}}
We will prove the conclusions of Lemma \ref{sigma} successively.

1. In view of the fact that the eigenvalues of $\bSigma$ and $\bPhi$ are the inverses of each other, the first conclusion is clear.

2. Suppose that the eigenvalues of $\bPhi$ are $1/L\leq \lambda_1\leq \lambda_2\leq \cdots\leq \lambda_p\leq L$. Then there is some orthonormal matrix $\bQ=(q_{ij})_{p\times p}$ such that
\begin{eqnarray}\label{toulan}
\bPhi=\bQ\left(\begin{array}{ccc}{\lambda_1}  & {\cdots} & {0} \\  {\vdots}  & {\ddots} & {\vdots} \\ {0} & {\cdots} & {\lambda_p}\end{array}\right)\bQ^\top.
\end{eqnarray}
Thus, for any $1 \leq k \leq p$, we have
$$1/L\leq \lambda_1\leq\phi_{kk}=\lambda_1q_{k1}^2+\ldots+\lambda_pq_{kp}^2\leq \lambda_p\leq L,$$
which completes the proof of this conclusion.

3. By the Cauchy's interlace theorem for eigenvalues in \cite[Theorem 1]{Hwang2004}, this conclusion is also clear.


4. Based on \eqref{toulan}, we have
\begin{eqnarray}\label{toulanb}
\bPhi\bPhi^\top=\bQ\left(\begin{array}{ccc}{\lambda_1^2}  & {\cdots} & {0} \\  {\vdots}  & {\ddots} & {\vdots} \\ {0} & {\cdots} & {\lambda_p^2}\end{array}\right)\bQ^\top.
\end{eqnarray}
Denote by $\|\bphi_j\|_2$ the $L_2$-norm of the $j$th row of $\bPhi$. Since $\bPhi$ is symmetric, $\|\bphi_j\|_2^2$ is equal to the $j$th diagonal component of $\bPhi\bPhi^\top$ for any $j\in \{1,\ldots,p\}$.

In view of \eqref{toulanb}, we know that the eigenvalues of the symmetric matrix  $\bPhi\bPhi^\top$ are bounded within the interval $[L^{-2}, L^2]$. Thus, the second conclusion of this lemma gives that $\|\bphi_j\|_2^2$ is bounded within the interval $[L^{-2}, L^2]$, which yields that $\|\bphi_j\|_2$ is bounded within the interval $[1/L, L]$. It completes the proof of this conclusion.

5. Based on \eqref{toulan}, we have
\begin{eqnarray*}
\bSigma=\bQ\left(\begin{array}{ccc}{1/\lambda_1}  & {\cdots} & {0} \\  {\vdots}  & {\ddots} & {\vdots} \\ {0} & {\cdots} & {1/\lambda_p}\end{array}\right)\bQ^\top.
\end{eqnarray*}
Thus, for any $1 \leq j \leq p$, it yields that
$$ \sigma_{jj}=\frac{1}{\lambda_1}q_{j1}^2+\ldots+\frac{1}{\lambda_p}q_{jp}^2.$$
By Cauchy-Schwartz inequality, we get
\begin{align*}
\sigma_{jj}\phi_{jj}&=\left(\frac{1}{\lambda_1}q_{j1}^2+\ldots+\frac{1}{\lambda_p}q_{jp}^2\right)(\lambda_1q_{j1}^2+\ldots+\lambda_pq_{jp}^2)\\
&=\left(\sum_{k=1}^p\left(\frac{q_{jk}}{\sqrt{\lambda_k}}\right)^2\right)\left(\sum_{k=1}^p (q_{jk}\sqrt{\lambda_k})^2\right)\geq \left(\sum_{k=1}^p q_{jk}^2\right)^2=1.
\end{align*}
It concludes the proof of Lemma \ref{sigma}.

\subsection{Proof of Lemma \ref{lse-lem}}
We first prove that with significant probability, $\bX^\top\bX$ is invertible so that the ordinary least squares estimator $\widehat{\bbeta}=(\bX^\top\bX)^{-1}\bX^\top\by$ exists. 
Clearly, $\bX\bSigma^{-1/2} \sim N(0, \bI_n \otimes \bI_d)$. Applying \cite[Theorem 5.39]{Eldar2012}, we can conclude that for some positive constants $c$ and $C$, with probability at least $1-2\exp(-cn)$, it holds that
\begin{align*}
\sqrt{n}-(C+1)\sqrt{d}\leq s_{\min}(\bX\bSigma^{-1/2})\leq s_{\max}(\bX\bSigma^{-1/2})\leq\sqrt{n}+(C+1)\sqrt{d},
\end{align*}
where $s_{\min}(\bX\bSigma^{-1/2})$ and $s_{\max}(\bX\bSigma^{-1/2})$ are the smallest and the largest singular values of $\bX\bSigma^{-1/2}$, respectively.
Since $d=o(n)$, it follows that for sufficiently large $n$, with probability at least $1-2\exp(-cn)$, $\bX\bSigma^{-1/2}$ will have full column rank.

Moreover, since $\bSigma$ is positive definite, we can deduce that with probability at least $1-2\exp(-cn)$, $\bX$ is of full column rank, so that $\bX^\top\bX$ is positive definite and thus invertible.
For any $p$ satisfying $\log(p)/n=o(1)$ and any constant $\delta > 1$, we have
$$\frac{\exp(-cn)}{p^{1-\delta}}=\frac{p^{\delta-1}}{\exp(cn)}=\exp\{(\delta-1)\log p - c n\}=o(1),$$
which entails $2\exp(-cn)=o(p^{1-\delta})$. Thus, we can claim that $\bX^\top\bX$ is positive definite and invertible with probability at least $1- o(p^{1-\delta})$.

We then show the prediction error bound of $\widehat{\bbeta}$. Since $\bX$ is independent of $\bveps$, by the properties of the ordinary least squares estimator, we have
\begin{eqnarray*}
\frac{\|\bX(\widehat{\bbeta}-\bbeta)\|_2^2}{\sigma^2}\big|\bX\sim \chi^2_{(d)}.
\end{eqnarray*}
Applying the tail probability bound of the chi-squared distribution in (\ref{chi}) with $t = 4\sqrt{\delta \log (p)/n}$ gives
$$\mathbb{P}\{d(1-t)\sigma^2   \le   \|\bX(\widehat{\bbeta}-\bbeta)\|_2^2 \le d(1+t)\sigma^2 \big|\bX \}\geq 1 - 2p^{-2\delta}.$$
Since the right hand side of the above inequality is independent of $\bX$, it follows that
$$
\mathbb{P}\{d(1+4\sqrt{\delta(\log p)/n})\sigma^2 \leq \|\bX(\widehat{\bbeta}-\bbeta)\|_2^2 \leq d(1+4\sqrt{\delta(\log p)/n})\sigma^2\}\ge 1- o(p^{1-\delta}),$$
which completes the proof of Lemma \ref{lse-lem}.

\subsection{Proof of Lemma \ref{lem2}}
Since $\mathcal{G}_{-}(\xi,S)\subset\mathcal{G}(\xi,S)$, it is immediate that $F_2(\xi,S)\geq F_2^* (\xi,S)$. To show the constant lower bound, we first establish a connection between $F_2 (\xi,S)$ and the restricted eigenvalues (RE), which are defined as
$$RE_m (\xi,S)= \inf \left\{\frac{s^{1/m}\|\bX \bu\|_{2}}{(ns)^{1/2} \|\bu\|_{m}}: \bu \in \mathcal{G}(\xi,S) \right\} $$ 
for a positive integer $m$ and the cone $\mathcal{G}(\xi,S)$ of the same definition as that in Lemma \ref{lem2}. By \cite[Proposition 5]{ye2010}, we have
$$F_2^* (\xi,S)\geq RE_1 (\xi,S)RE_2 (\xi,S).$$
Since $RE_1 (\xi,S)\geq RE_2 (\xi,S)$, it follows that
$$F_2^* (\xi,S)\geq RE_2^2 (\xi,S).$$

We then prove a constant lower bound for $RE_2 (\xi,S)$ with significant probability. First, it follows from \cite[Theorem 1]{Raskutti2010} that for any $\bu\in \mathbb{R}^p$, there exist positive constants $c_2$ and $c_3$ so that
\begin{align*}
\frac{\|\bX \bu\|_{2}}{n^{1/2}}\geq \frac{1}{4}\|\mathbf{\Sigma}^{1/2}\bu\|_2-9\rho(\mathbf{\Sigma})\sqrt{\frac{\log p}{n}}\|\bu\|_1 
\end{align*}
holds with probability at least $1- c_2\exp(-c_3n)$, where $\rho(\mathbf{\Sigma})=\max_{j\in \{1,\ldots,p\}}\mathbf{\Sigma}_{jj}^{1/2}$. 
On the other hand, for $\bu\in \mathcal{G}(\xi,S)$, it follows from the definition of $\mathcal{G}(\xi,S)$ that
$$\|\bu\|_1=\|\bu_{S^c}\|_1 + \|\bu_{S}\|_1\leq (1+\xi) \|\bu_{S}\|_1\leq (1+\xi)s^{1/2} \|\bu_{S}\|_2.$$

In addition, we have $\|\mathbf{\Sigma}^{1/2}\bu\|_2\geq \lambda_{\min}(\mathbf{\Sigma}^{1/2})\|\bu\|_2$, where $\lambda_{\min}(\mathbf{\Sigma}^{1/2})$ is the smallest eigenvalue of $\mathbf{\Sigma}^{1/2}$. Thus, it follows that
\begin{align}\label{wuyule}
\frac{\|\bX \bu\|_{2}}{n^{1/2}}\geq \left\{\frac{1}{4}\lambda_{\min}(\mathbf{\Sigma}^{1/2})-9(1+\xi)\rho(\mathbf{\Sigma})\sqrt{\frac{s \log (p)}{n}}\right\}\|\bu\|_2
\end{align}
holds with probability at least $1- c_2\exp(-c_3n)$. Since the eigenvalues of $\bSigma$ are bounded within the interval $[1/L, L]$, both $\lambda_{\min}(\mathbf{\Sigma}^{1/2})$ and $\rho(\mathbf{\Sigma})$ are bounded from above and below by some positive constants according to Lemma \ref{sigma}.

Therefore, when $s=o(n/\log (p))$, it follows from the definition of $RE_2 (\xi,S)$ that there exists some positive constant $c_1$ such that for sufficiently large $n$,
\begin{align*}
RE_2 (\xi,S)  \geq c_1^{1/2}
\end{align*}
holds with probability at least $1- c_2\exp(-c_3n)$. It completes the proof of Lemma \ref{lem2}.

\end{document}